\documentclass[english,reqno]{amsart}

\usepackage{amsmath}
\usepackage{amssymb} %to produce blackboard bold characters(\supsetneq)
\usepackage{enumitem} %lay­out of the three ba­sic list en­vi­ron­ments
\usepackage{mathtools} %an extension package to amsmath + loading amsmath
\usepackage[table]{xcolor} %color of tables
\usepackage[all]{xy} %xymatrix commutative diagram
\usepackage{tikz} %pictures
\usepackage{tikz-cd}%implies
\usepackage{indentfirst} %Indent first paragraph
\usepackage{babel} %lan­guages
\usepackage{setspace} %set­ting the spac­ing be­tween lines
\usepackage[colorlinks,linkcolor=red,anchorcolor=green,citecolor=blue]{hyperref} %color of hyperlinks
\hypersetup{linktocpage = true} %color of hyperlinks of TOC

\usepackage{rotating} %Rotate an image, table or paragraph
\usepackage{ytableau} %young tableau
\usepackage{longtable} %long table
\newcolumntype{M}[1]{>{\centering\arraybackslash}m{#1}} %define dimension for long stable

\usepackage{mathpazo}
\usepackage{mathrsfs}
\DeclareFontFamily{OMS}{rsfs}{\skewchar\font'60}
\DeclareFontShape{OMS}{rsfs}{m}{n}{<-5>rsfs5 <5-7>rsfs7 <7->rsfs10 }{}
\DeclareSymbolFont{rsfs}{OMS}{rsfs}{m}{n}
\DeclareSymbolFontAlphabet{\scr}{rsfs}
\DeclareSymbolFontAlphabet{\scr}{rsfs}

\usepackage[T1]{fontenc}
\theoremstyle{plain}
\newtheorem{thm}{Theorem}[section]
\newtheorem{lemma}[thm]{Lemma}
\newtheorem{prop}[thm]{Proposition}

\newtheorem{defn}[thm]{Definition}

\theoremstyle{remark}

\newtheorem{remark}[thm]{Remark}

\def\Im{\operatorname{Im}}

\def\div{\operatorname{div}}
\def\mod{\operatorname{mod}}

\def\dim{\operatorname{dim}}

\def\Proj{\operatorname{Proj}}

\def\ker{\operatorname{ker}}

\def\min{\operatorname{min}}
\def\max{\operatorname{max}}

\def\hor{\operatorname{\hor}}
\def\ver{\operatorname{\ver}}

\def\sm{\operatorname{\textsubscript{\rm sm}}}
\def\sing{\operatorname{\textsubscript{\rm sing}}}

\def\ver{\operatorname{\textsubscript{\rm ver}}}
\def\hor{\operatorname{\textsubscript{\rm hor}}}
\def\rank{\operatorname{rank}}

\def\ddc{\frac{\sqrt{-1}}{\pi} \partial \bar{\partial}}

\setlist[itemize]{leftmargin=*}
\setlist[enumerate]{leftmargin=*}

\numberwithin{equation}{section} %numbering of equations
\makeatletter

\ifnum\@ptsize=0  \fi \ifnum\@ptsize=2  \fi

\title{Title} 
\subjclass[2010]{}

\author{Wenhao Ou}

\address{Wenhao Ou, Institute of Mathematics, Academy of Mathematics and Systems Science, Chinese Academy of Sciences, Beijing, 100190, China}
\email{wenhaoou@amss.ac.cn}

\begin{document}

\begin{abstract}
We adapt Bost's algebraicity characterization to the situation of a  germ in a compact K\"ahler manifold. 
As a consequence, we extend the  algebraic integrability criteria   of Campana-P{\u a}un and of Druel to foliations on compact K\"ahler manifolds. 
As an application, we prove that a compact K\"ahler manifold  is uniruled  if and only if its canonical line bundle is not pseudoeffective.
\end{abstract}

\title{A characterization of uniruled compact K\"ahler manifolds}

\maketitle

\tableofcontents

%\vspace{-0.2cm} 

\section{Introduction}

A compact K\"ahler manifold $X$ is called uniruled if it is dominated by a family of rational curves. 
Uniruled varieties form an important class in the classification of varieties. 
It is conjectured that  $X$ is uniruled if and only if its Kodaira dimension is negative. 
In \cite{MiyaokaMori1986}, Miyaoka and Mori showed that a projective manifold $X$ is uniruled if and only if it is dominated by a family of curves $C$, such that the intersection numbers $\omega_X \cdot C$ are negative, 
where $\omega_X$ is the canonical line bundle of $X$.  
In \cite{BoucksomDemaillyPuaunPeternell2013}, 
Boucksom, Demailly, P{\u a}un and Peternell showed that the latter condition is equivalent to that $\omega_X$ is not pseudoeffective.  
It is natural to expect an analogue characterization holds for compact K\"ahler manifolds, \textit{e.g.} \cite{CaoHoering2020}, \cite{HaconPaun2024}.  
In this paper, we prove the following theorem which confirms the conjecture. 
We also refer to \cite{Yau1974} for the case of surfaces,
and \cite{Brunella2016} for the case of threefolds.

\begin{thm}
\label{thm:uniruledness}
Let $X$ be a compact K\"ahler manifold. 
Then $X$ is uniruled if and only if the canonical line bundle $\omega_X$ is not pseudoeffective. 
\end{thm}

The method of \cite{Brunella2016} proceeds as  follows. 
If $X$ is a non projective compact K\"ahler threefold, 
then there is a non zero holomorphic $2$-form  on $X$, 
whose kernel induces a foliation $\mathcal{F}$ of rank $1$. 
If $\omega_X$ is not pseudoeffective, 
then $\mathcal{F}^*$ is not pseudoeffective neither. 
Brunella managed to show that the closures of the leaves of $\mathcal{F}$ are rational curves. 
Our proof of  Theorem \ref{thm:uniruledness}  follows the same strategy, 
and the crux is to prove some criteria for foliations induced by meromorphic maps (which are algebraically integrable foliations in the algebraic setting).  
It is  related to the following problem.

Let $X$ be a compact complex manifold and let $S_0\subseteq X$ be  an irreducible locally closed submanifold. 
It is  natural to investigate if $S_0$ can extend to a closed complex analytic subvariety $Y$ of $X$, so that $S_0$ is a local branch of $Y$.  
This is equivalent to that the Zariski closure of $S_0$ has the same dimension as $S_0$.  
For example, we assume that  $S_0\subseteq \mathbb{C}^2$ is a Euclidean open subset of the graph $\{y=f(x)\}$, where  $f$  is some holomorphic function on $\mathbb{C}$.    
Then  the Zariski closure  of $S_0$ in   $X=\mathbb{P}^2$ is a curve if and only if $f$ is a polynomial,  
for   closed analytic subvarieties of $\mathbb{P}^2$ are algebraic by Chow's theorem.  

 %$S_0$ is the analytic graph of a foliation $\mathcal{F}$ on some projective manifold,  and if $S_0$ is  has the same dimension as its Zariski closure,  then $\mathcal{F}$ is algebraically integrable. 

%If it is the case, we  say that $S_0$ is openin its Zariski closure.  
When $X$ is projective, Bost proved a characterization for the algebraicity of a germ $(x\in S_0)$ in \cite{Bost04} as follows.    
Let $\mathcal{L}$ be an ample line bundle on $X$. 
For any integers $D, i > 0$, we define the vector subspace $E_D^i$ of $H^0(X, \mathcal{L}^{\otimes D})$ as the collection of sections whose restrictions on $S_0$ vanish with orders at least $i$ at $x$.  
Then $S_0$ has the same dimension as its Zariski closure if and only if there is some constant $\lambda >0$ such that $E_{D}^i=E_{D}^{i+1}$ whenever  $i \cdot D^{-1}>  \lambda$. 

We observe that, in Bost's method, a section $\sigma$ of $\mathcal{L}^{\otimes D}$ induces canonically a singular Hermitian metric $h$ on $\mathcal{L}$. 
More precisely, if  $\omega$ is a K\"ahler form with class in $c_1(\mathcal{L})$, 
then there is  a $\omega$-psh function $\varphi$  of the shape $\varphi = \frac{1}{D} \log |\sigma| + O(1)$.  
If $\sigma \in E_{D}^i \setminus E_{D}^{i+1}$,   
then the term $i\cdot D^{-1}$ is equal to the Lelong number $\nu(\varphi|_{S_0}, x)$. 
With this perspective, one may expect to extend Bost's characterization to compact K\"ahler manifolds, in terms of quasi-psh functions and of Lelong numbers. 
We first prove the following theorem.

\begin{thm}
\label{thm:Zariski-dense-germ}
Let $(X,\omega)$ be a compact K\"ahler manifold, 
let $C\subseteq X$ be a closed irreducible submanifold, 
and let $S_0\subseteq X$ be an irreducible locally closed subvariety containing a Zariski open subset $C_0\neq \emptyset$ of $C$.   
We make the following assumptions. 
\begin{enumerate}
    \item[(i)]   $S_0$ is Zariski dense in $X$.  
    \item[(ii)]  $\dim S_0 = \dim C + 1$ %and the normal bundle of $C_0$ in $S$ extends to a  line bundle  $\mathcal{N}$ on $C$.  
    \item[(iii)]  The prime divisors in $C$   contained in $C\setminus {C_0}$  form an exceptional family.  
    \item [(iv)] %$(S_0,C)$ admits infinitely many blowups. 
       There is a Zariski open subset $C_1$ of $C_0$ whose complement has codimension at least $2$, such that $S_0$ is smooth around $C_1$. 
       Furthermore, $S_0$ extends formally along $C$. 
\end{enumerate}
Then for any  $\lambda>0$, there is a $\omega$-psh function $\varphi$ satisfying the following properties. 
\begin{enumerate} 
    \item $\varphi$ has analytic singularities. In addition,   $\varphi$ can be locally written as 
    \[\varphi = \frac{1}{2m} \log  (|g_1|^2 +  \cdots + |g_r|^{2}) + O(1)\] 
    where $m\in \mathbb{Z}_{>0}$,
    and  $g_1,...,g_r$ are local  generators of  a coherent ideal sheaf $\mathcal{J}$ on $X$. 
    \item For all points $y \in C_1$, we have $\nu(\varphi|_{S_0}, y) > \lambda$.
\end{enumerate}
\end{thm}

We refer to Section \ref{section:Zariski} for the notion of  exceptional family in the assumption (iii), and Section  \ref{section:blowup} for the notion   in the assumption (iv). 
Particularly, these assumptions are satisfied if $C_0=C$.   
With Theorem \ref{thm:Zariski-dense-germ}, we can prove the following theorem, 
which is an extension of the algebraicity criteria in \cite{BogomolovMcQuillan16} and in \cite{Bost01}.

\begin{thm}
\label{thm:alg-germ-nonpsef-normal} 
Let $X$ be a compact K\"ahler manifold, 
let $C\subseteq X$ be a closed irreducible submanifold, 
and let $S_0\subseteq X$ be an irreducible locally closed submanifold. 
Assume the following conditions. 
\begin{enumerate}
\item $S_0$ contains a   Zariski open subset $C_0\neq \emptyset$ of $C$. 
\item %$(S_0,C)$ admits infinitely many blowups.  
       $S_0$ extends formally along $C$. 
      In particular, the conormal bundle  $\mathcal{N}_{C_0/S_0}^*$ extends to a reflexive coherent sheaf $\mathcal{N}^*$ on $C$.  
\end{enumerate}
In addition, we assume one of the following conditions holds, 
\begin{enumerate}
\item[(3)] either $\dim S_0 = \dim C+1$,   
the prime divisors in $C$ contained in $C\setminus C_0$  form an exceptional family,  
and for any divisor  class $\delta $ supported in $C\setminus C_0$,  
the class $c_1( \mathcal{N}^{*}) + \delta$ is not pseudoeffective; 
\item[(4)]  or  $ C\setminus C_0$ has codimension at least 2, 
and $\mathcal{N}^*$ is non pseudoeffective in the sense of Definition \ref{def-non-psef}.   
\end{enumerate}
Then $S_0$ is has the same dimension as its Zariski closure. 
\end{thm}

In the  simple  case when $X$ is a surface, $C=\{x\}$ is a point, 
and $S_0$ is a germ of curve around $x$, 
the proof of Theorem \ref{thm:Zariski-dense-germ} is sketched as follows.   
By Demailly's    mass concentration in \cite{Demailly1993}, 
we can construct a $\omega$-psh function $\psi$ with arbitrarily large Lelong number $\nu(\psi|_{S_0}, x)$. 
Then we apply Demailly's regularization of current in \cite{Demailly1992} to $\psi$ and obtain some 
$\omega$-psh function $\varphi$.   
The construction of the regularization then implies that $\varphi$ satisfies the desired property.

Let $\mathcal{F}\subseteq T_X$ be a foliation on a compact K\"ahler manifold $X$. 
We say that $\mathcal{F}$ is induced by a meromorphic map  if it is induced by a dominant meromorphic map $X\dashrightarrow Y$,  where $Y$ is a  compact complex analytic variety. 
Assume that  $X_0$ is the largest open subset where $\mathcal{F}$ is a subbundle. 
Let $\Delta \subseteq X\times X$ be the diagonal  and $\Delta_0=\Delta \cap (X_0\times X_0)$. 
The analytic graph of $\mathcal{F}$ is a locally closed submanifold ${S_0}\subseteq X_0\times X_0$ containing $\Delta_0$.  
Then the foliation $\mathcal{F}$ is  induced by a meromorphic map 
if   ${S_0}$ has the same dimension as its Zariski closure in $X\times X$, 
see Lemma \ref{lemma:foliation-graph}.  
By using Theorem \ref{thm:alg-germ-nonpsef-normal}, 
we prove the following theorem, 
which is  the key ingredient for Theorem \ref{thm:uniruledness}. 
It extends the algebraic  integrability criteria of \cite{CampanaPaun2017} and of \cite{Druel2018} (we also refer to \cite{Miyaoka1987}, \cite{BogomolovMcQuillan16}, \cite{Bost01},  \cite{Bost04} and \cite{Brunella2016} for more classical results).  
In our context,  
a class $\alpha\in H^{n-1,n-1}(X,\mathbb{R})$ is called movable if and only if it has non negative intersections with pseudoeffective classes. 

\begin{thm}
\label{thm:foliation-psef-alg} 
Let $X$ be a compact K\"ahler manifold of dimension $n$,
and let $\mathcal{F}$ be a foliation on $X$. 
\begin{enumerate}
    \item If  $\mathcal{F}^{*}$ is non pseudoeffective in the sense of Definition \ref{def-non-psef}, then $\mathcal{F}$ is  induced by a meromorphic map. 
    \item In particular, if the minimal slope satisfies  $\mu_{\alpha,min}(\mathcal{F}) >0$  
    for some movable class $\alpha \in H^{n-1,n-1}(X,\mathbb{R})$, then $\mathcal{F}$ is  induced by a meromorphic map. 
\end{enumerate} 
\end{thm}

We would like to thank St\'ephane Druel for pointing out the following application of Theorem \ref{thm:alg-germ-nonpsef-normal}, 
which extends  the algebraicity criterion of  
\cite{BogomolovMcQuillan16} (see also \cite{KebekusSolaToma}).
%\cite[Main Theorem]{BogomolovMcQuillan16} (see also \cite[Theorem 1]{KebekusSolaToma}). 

\begin{thm}
\label{thm:BMQ} 
Let $X$ be a compact K\"ahler variety,  let $\mathcal{F}$ be a foliation on $X$,  
and let $V\subseteq X$ be an irreducible subvariety contained in the regular locus of $\mathcal{F}$. 
Assume that, for some desingularization $\rho\colon W \to V$, the vector bundle $ \rho^*\mathcal{F} $ is non pseudoeffective in the sense of Definition \ref{def-non-psef}.  
Then for any point $x\in V$,  
the local leaf of $\mathcal{F}$ passing through $x$ has the same dimension as  its Zariski closure in $X$. 
\end{thm}

\vspace{5mm} 

{\bf Acknowledgments.} 
The author would like to thank Jian Xiao for informing him of the paper \cite{Demailly1993}.
He is grateful to 
Junyan Cao, 
St\'ephane Druel,  
Christopher Hacon 
and Burt Totaro 
for reading preliminary versions of the paper, 
and providing a lot of valuable advice.  
The author is supported by the National Key R\&D Program of China (No. 2021YFA1002300).

\section{Preliminaries}  

Throughout this paper, 
a  {complex analytic variety}   $X$
is a reduced and irreducible  complex space.   
We will work with Euclidean topology unless otherwise specified. 
By a  general point of $X$, 
we refer to a point which belongs to some non empty Zariski open subset of $X$. 
We say that a general point of $X$ satisfies certain property $\mathscr{P}$, 
if there is some non empty Zariski open subset of $X$, such that each point of it satisfies the property $\mathscr{P}$.  
We denote the  smooth locus of $X$  by $X_{\sm}$ and the singular locus by $X_{\sing}$.  
A complex manifold  is a smooth complex analytic variety.    
We will call a locally free coherent sheaf a vector bundle. 
We denote by $\Omega^1_X$   the sheaf of K\"ahler differentials on $X$, 
by $T_X = (\Omega_X^1)^{*}$ the reflexive tangent sheaf, 
and by $\omega_X= (\bigwedge^{\dim X} \Omega_X^1)^{**}$ the canonical sheaf.   
If $X$ is smooth, then $\Omega_X^1$ is  equal to the (holomorphic) cotangent bundle of $X$. 
In this case,  we set  $\Omega_X^r = \bigwedge^r \Omega_X^1$ for any integer $r\ge 1$.

A morphism $f\colon X\to Y$ between complex analytic varieties is called a fibration if it is proper surjective  and if $f_*\mathcal{O}_X \cong \mathcal{O}_Y$. 
In particular, every fiber of a fibration is connected. 
We denote by $\omega_{X/Y}= (\omega_X\otimes (f^*\omega_Y)^*)^{**}$ the relative canonical sheaf.

%We always work with complex analytic varieties $X$ which are homotopy equivalent to some finite CW-complexes, 
%so that   their singular  cohomology groups with coefficients in $\mathbb{Z}$ are finitely generated. 
Let $X$ be a complex manifold. 
We denote by $H^i(X,A)$   the singular cohomology groups with coefficients in an Abelian group $A$. 
If $A$ is a field of characteristic zero, 
and if $H^i(X,\mathbb{Z})$ is finitely generated, 
then the universal coefficient theorem implies that the natural morphism 
\[
H^i(X,\mathbb{Z}) \otimes_{\mathbb{Z}} A \to  H^i(X,A)  
\]
is isomorphic. 
If $A=\mathbb{R}$ or $\mathbb{C}$, then de Rham's theorem states that the de Rham cohomology groups are canonically isomorphic to the singular cohomology groups.  
Therefore, by abuse of notation, we also write $H^i(X,\mathbb{R})$ (respectively  $H^i(X,\mathbb{C})$) for the de Rham cohomology groups of  real (respectively complex) differential forms or currents.  
For a closed differential form or a closed current $\alpha$, 
we write $\{\alpha\} \in H^i(X,\mathbb{Q})$  (respectively  $\{\alpha\} \in H^i(X,\mathbb{Z})$), 
if the  class $\{\alpha\}$ belongs to the  image of $H^i(X,\mathbb{Q})$ (respectively  of $H^i(X,\mathbb{Z})$).  
If   $X$ is a compact K\"ahler manifold, then there is a canonical Hodge decomposition 
\[
H^i(X,\mathbb{C}) = \bigoplus_{p+q = i} H^{p,q}(X,\mathbb{C}), 
\]
where $H^{p,q}(X,\mathbb{C}) \subseteq H^i(X,\mathbb{C})$ is the subspace of classes represented by closed $(p,q)$-forms.  
Moreover, there is a canonical isomorphism  $H^{p,q}(X,\mathbb{C}) \cong H^q(X,\Omega_X^p)$ for all $p,q\ge 0$. 
%In addition, there is a natural isomorphism $H^{p,q}(X,\mathbb{C}) \cong \overline{H^{p,q}(X,\mathbb{C})}$, where the overline signifies the complex conjugate.   
For any integer $j\ge 0$, we denote by $H^{j,j}(X,\mathbb{R})$ the intersection $ H^{j,j}(X,\mathbb{C})\cap H^{2j}(X,\mathbb{R})$.

\begin{lemma}
\label{lemma:torsion-line-bundle}
Let $X$ be a complex manifold and let $\mathcal{L}$ be a line bundle on $X$.  
We assume that $H^2(X,\mathbb{Z})$ is finitely generated,  
that $H^1(X,\mathcal{O}_X) = \{0\}$, 
and that $c_1(\mathcal{L})=0 \in H^{2}(X,\mathbb{R})$.  
Then there is an integer $m>0$ such that $\mathcal{L}^{\otimes m} \cong \mathcal{O}_X$. 
\end{lemma}

\begin{proof}
The standard exponential exact sequence 
\[
0\to \mathbb{Z} \to \mathcal{O}_X \to \mathcal{O}_X^\times \to 0
\]
induces an exact sequence 
\[
H^1(X,\mathcal{O}_X) \to H^1(X,\mathcal{O}_X^\times)  \overset{\iota}{\to} H^2(X,\mathbb{Z}). 
\]
The middle term above is the Picard group of $X$, 
and the composition of $\iota$ with the natural morphism $p\colon H^2(X,\mathbb{Z}) \to H^2(X,\mathbb{R})$ is equal to taking the first Chern class. 
The assumption implies that  $\iota$ is an injective morphism. 
Hence  $\ker p\circ \iota$ is isomorphic to a subgroup of $ \ker p$, 
and is an Abelian torsion group. 
Therefore, there is some $m>0$ such that such that $\mathcal{L}^{\otimes m} \cong \mathcal{O}_X$.  
\end{proof}

\begin{lemma}
\label{lemma:extend-class} 
Let $X$ be a compact complex  manifold, 
let $Z\subseteq X$ be a closed analytic  subset,  
let $X_0 = X\setminus Z$ and let $U\subseteq X_0$ be a Zariski open subset. 
We suppose that the codimension of $X_0\setminus U$ in $X_0$ is at least $2$.  
Assume that there are two line bundle $\mathcal{L}$ and $\mathcal{L}'$ on $X$ such that 
$\mathcal{L}'|_U \cong \mathcal{L}|_U$. 
Then there is a divisor class $\delta$ with support contained in $Z$, 
such that 
\[
c_1(\mathcal{L}) = c_1(\mathcal{L}')+\delta  \in H^2(X,\mathbb{R}). 
\]
\end{lemma} 

We remark that $U$ may not be Zariski open in $X$. 

\begin{proof} 
Without loss of the generality, we can assume that $\mathcal{L}'=\mathcal{O}_X$. 
By assumption, there is a nowhere vanishing section $\sigma\in H^0(U,\mathcal{L}|_U)$. 
Since $X$ is smooth and since  the codimension of $X_0\setminus U$ in $X_0$ is at least $2$,  by working locally at around every point of $X_0$, 
we deduce that $\sigma$ extends to a nowhere zero section of $\mathcal{L}|_{X_0}$. 
This implies that  $\mathcal{O}_{X_0}\cong \mathcal{L}|_{X_0}$.    
By using the same argument, we may assume further that 
$Z$ is not empty and  has pure codimension $1$.

The previous paragraph implies that $c_1(\mathcal{L})$ lies in the image of the relative cohomology group $H^2(X,X_0,\mathbb{R})$.  
Let $n=\dim X$.  Since $X$ is smooth and compact, the group $H^2(X,X_0,\mathbb{R})$ is isomorphic to the homology group $H_{2n-2}(Z,\mathbb{R})$ by taking the cap product with the fundamental class of $X$. 
Since $Z$ has pure dimension $n-1$, 
we see that $H_{2n-2}(Z,\mathbb{R})$ is generated by the fundamental classes of the irreducible component of $Z$.
Thus $c_1(\mathcal{L}) = \delta$ for  a divisor class $\delta$ with support contained in $Z$. 
This completes the proof of the lemma. 
\end{proof}

A foliation $\mathcal{F}$ on a complex analytic variety $X$ is a saturated coherent subsheaf of $T_X$, which is closed under the Lie bracket. 
The regular locus of $\mathcal{F}$ is the largest open subset of $X_{\sm}$ where $\mathcal{F}$ is a subbundle of $T_X$. 
The singular locus of $\mathcal{F}$, denote by $\mathcal{F}_{\sing}$, is defined as the complement of the regular locus of $\mathcal{F}$.  
We remark that $\mathcal{F}_{\sing}$ is a closed analytic subset of $X$.

Assume that $\theta$ is a smooth real $(1,1)$-form on a  complex manifold $X$. 
A function $\varphi\colon X\to [-\infty, +\infty)$ is called a $\theta$-psh function if it is upper semicontinuous, locally integrable  
and  if $\theta + \ddc\varphi $ is  positive in the sense of currents. 
The word psh is short for plurisubharmonic. 
A function is called quasi-psh or almost psh if it is $\theta$-psh for some smooth  real $(1,1)$-form $\theta$.   
The Lelong number of a quasi-psh function $\varphi$ at a point $x\in X$ is defined as 
\[
\nu(\varphi, x) = \sup \{ \lambda \ge 0 \ | \  \varphi (y) \le \lambda \log |y-x| + O(1) \mbox{ around } x \}, 
\]
where $|y-x|$ is the Euclidean distance between $y$ and $x$ on some coordinate chart around $x$.  
The multiplier ideal sheaf $\mathcal{I}(\varphi)$ of $\varphi$ is defined as the coherent ideal sheaf locally generated by the holomorphic functions $f$ such that $|f|^2e^{-2\varphi}$ are locally integrable.

If $X$ is compact, then every positive closed $(1,1)$-current $T$ can be written as $\theta + \ddc \varphi$ where $\theta$ is a smooth form and $\varphi$ is a $\theta$-psh function. 
In this case,   the Lelong number of $T$ at $x$ is equal to the one of $\varphi$ at $x$.  
We recall Siu's decomposition of positive currents.  

\begin{thm}
\label{thm-decomposition}
Let $X$ be a   complex manifold and let $T$ be a positive  closed $(1,1)$-current. 
Then $T = R + \sum_{i\in I} a_i [E_i],$ for some index set  $I$. 
The $E_i$  are pairwise distinct prime divisors in $X$, 
$0<a_i=\inf_{x\in E_i}  {\nu}(T, x)$, 
and $R$ is a positive current such that $\inf_{x\in D}  {\nu}(T, x)=0$ for any prime divisor $D$. 
\end{thm}

\begin{proof}
See \cite{Siu1974} or \cite[(2.18)]{Demaillybook2012}. 
\end{proof}

A quasi-psh function $\varphi$ on $X$ is said to have analytic singularities around a point $x$ if locally around $x$, 
we can write 
\[
\varphi = \frac{\alpha}{2} \cdot \log (|g_1|^2+\cdots |g_r|^2) + O(1)
\]
where $\alpha>0$ is a real number and  $g_1,...,g_r$ are holomorphic functions around $x$.   
In this case,   $\nu(\varphi, x)$ is equal to $\alpha$ multiplied by the minimal vanishing order of the $g_i$ at $x$ .  
A closed positive   $(1,1)$-current $T$ is said to have analytic singularities around $x$ if it can be locally written as $\theta+\ddc \varphi$ with $\theta$ smooth such that  $\varphi$ has analytic singularities. 
The following statement follows from Demailly's regularization theorem.  

\begin{thm}
\label{thm:regularization}
Let $(X,\omega)$ be a compact K\"ahler manifold of dimension $n$, 
let $\theta$ be a smooth $(1,1)$-form and let $\psi$ be a $\theta$-psh function. 
Then for every $\varepsilon >0$, for any integer $m_0>0$, 
there is a quasi-psh function $\varphi$ satisfying the following properties. 
\begin{enumerate}
    \item $\theta+\ddc \varphi + \varepsilon \omega \ge 0$. 
    \item There is an integer $m>m_0$, such that $\varphi$ can be locally written as 
    \[
    \varphi = \frac{1}{2m} \log (|g_1|^2+\cdots+|g_r|^2) + O(1),
    \]
    where $g_1,...,g_r$ are local generators of a coherent ideal sheaf $\mathcal{J}$ on $X$. 
    \item There is an open covering $\{U_i\}$ of $X$, such that 
    $\mathcal{J}$ on each $U_i$ is generated by a basis of the Hilbert space 
    $\mathscr{H}_{U_i}(m\psi_i) = \{f\in \mathcal{O}_X(U_i) \  | \  \int_{U_i} |f|^2e^{-2m\psi_i} < +\infty \}$, 
    where $\psi_i = \psi|_{U_i}-q_i$ for some smooth function $q_i$ on $U_i$.  
    In particular, $\mathcal{J}\subseteq \mathcal{I}(m\psi)$. 
    \item $\nu(\psi,x) - \frac{n}{m} \le \nu(\varphi, x) \le \nu(\psi, x)$ for all $x\in X$. 
\end{enumerate}
\end{thm}

\begin{proof}
The construction follows from \cite[Proposition 3.1]{Demailly1992}.  
The existence of  the global ideal sheaf $\mathcal{J}$  is discussed in the paragraph after \cite[Corollary 13.13]{Demaillybook2012}. 
\end{proof}

When $T$ is a positive current with analytic singularities around a point $x$ in a complex manifold $X$, we can ``extract'' the Lelong number $\nu(T,x)$ by blowing up $x$. More precisely, if $\rho \colon X'\to X$ is the blowup of $x$, then  $R:=\rho^*T-\nu(T,x) [E]$ is a positive current, 
and  a local potential of it is bounded around a general point of $E$.  
In the following lemma, we study the situation of $\nu(T|_S, x)$, where $S\subseteq X$ is a locally closed submanifold containing $x$. 
We remark that  $\nu(T|_S, x)$ can be much larger than $\nu(T,x)$. 

\begin{lemma}
\label{lemma:blowup-lelong} 
Let $ C\subseteq S \subseteq X$ be germs of complex manifolds around a point $o$. 
Assume that $\dim S = \dim C + 1$. 
Let $\varphi = \frac{1}{2} \log ( |g_1|^2 + \cdots + |g_r|^2)$ such that $g_1,...,g_r$ are holomorphic functions on $X$. 
Assume that %all $g_i$ vanish  along $C$, and 
there is at least one $g_i$ which is non zero at general points of $S$.  
Suppose that $\nu(\varphi|_S, x) = l $ for a general point $x\in C$. 
Then there is a projective bimeromorphic morphism $\rho \colon X'\to X$ satisfying the following properties. 
\begin{enumerate}
    \item $\rho$ is obtained by firstly blowing up $C$,  then blowing up the intersection of the strict transform of $S$ and the exceptional locus for several times. 
    \item The strict transform $S'$ of $S$ in $X'$ meets the $\rho$-exceptional divisor transversally at a unique irreducible component $E$. If $C'$ is the intersection of $S'$ and $E$, and if $o'\in C'$ is the  point lying over $o\in C$, 
    then     \[\varphi \circ \rho = {l} \cdot log |h| + \frac{1}{2} \log v   \] in a neighborhood of $o'$, 
    where $h$ is a local generator of   $\mathcal{O}_{X'}(-E)$, 
    and $v$ is a finite sum of squares of absolute values of   holomorphic functions. 
    Moreover, $v$  is not identically zero along $C'$.  
\end{enumerate} 
\end{lemma}

\begin{proof}
%If $l=0$, then we can let $\rho$ be the blowup of $X$ at $C$. 
%From now on we assume that $l>0$. 
%Then all $g_i$ vanish  along $C$. 
We may assume that $X$ is an open neighborhood of the origin $ o$ in $ \mathbb{C}^n$, and that $C$ and $S$ are coordinates subspaces. 
More precisely, we can assume that $(z_1,...,z_k,t,u_{k+2},...,u_n)$ is a coordinates system  of $\mathbb{C}^n$ and that 
\[
C=\{z_1=\cdots=z_k=t=0\},  \ \ \ S= \{z_1= \cdots z_k=0\}. 
\] 
We write $\mathbf{z}=(z_1,...,z_k)$ and $\mathbf{u}=(u_{k+2},...,u_n)$. 
%Since every $g_i$ vanishes along $C$, 
Up to shrinking $X$ around $o$, 
we can assume that each $g_i$ admits a decomposition as follows, 
\[
g_i(\mathbf{z},t,\mathbf{u}) 
= h_i(\mathbf{z},t,\mathbf{u}) + t^{a_i}\cdot \eta_i(\mathbf{z},t,\mathbf{u}),  
\]
where $h_i$ is a holomorphic function vanishing along $S$, $a_i \ge 0$ is an integer and $\eta_i$ is either   zero 
or a holomorphic function which is not identically zero along $C$. 
Since there is at least one $g_i$ which does not vanish entirely along $S$, 
there is at least one $\eta_i$ which is non zero.  
It follows that, for a general point $x\in C$,   
\[
l= \nu(\varphi|_S, x) = \min \{a_i \ | \ \eta_i \neq 0\}. 
\]
Without loss of the generality, we can assume that $a_i \ge l$ for all $i=1,...,r$.

Assume that $\rho\colon X'\to X$ is a composition of $m$ blowups, as described in the item (1) of the lemma.   
Then there is an open neighborhood $U'$ of $o'$, such that $\rho$ can be written in coordinates as 
\[
(z_1',...,z_k',t',u_{k+2}',...,u_n') \mapsto  
(t'^m z_1',...,t'^m z_k',t',u_{k+2}',...,u_n').  
\]
In addition, $\mathcal{O}_{X'}(-E)$ is generated by $t'$ on $U$, and 
\[
C'=\{z_1'=\cdots=z_k'=t'=0\},  \ \ \ S'= \{z_1'= \cdots z_k'=0\}. 
\] 
Since $h_i$ vanishes along $S$, as a series, every term of it is divisible by some $z_j$, with $j=1,...,k$. 
Hence there is some holomorphic function $h'_i$ on $U'$ such that 
\[
g_i\circ \rho (\mathbf{z}',t,\mathbf{u}') 
= t'^m h'_i(\mathbf{z}',t,\mathbf{u}') + t'^{a_i}\cdot \eta_i'(\mathbf{z}',t',\mathbf{u}'), 
\]
where $\eta_i'=\eta_i\circ \rho$.  

In particular, if we let $m=l+1$, then  $ g_i\circ \rho  
= t'^l ( t'\cdot h'_i  + t'^{a_i-l}\cdot \eta_i')$. 
By assumption, there is some $e$ such that $\eta_{e}'$ is non zero at  a  general point  $x'$  of $C'\cap U$ and that  $a_{e}=l$.  
Therefore, $|t'\cdot h'_{e} + t'^{a_e -l }\eta'_{e}|^2 >0$ at $x'$. 
Let 
\[v=  \sum_{i=1}^r | t'\cdot h'_i + t'^{a_i-l}\cdot \eta'|^2.\]
Then  $ \varphi \circ \rho = l\log |t'| +  \frac{1}{2}\log v$ in a neighborhood of $x'$.  
This completes the proof of the lemma.  
\end{proof}

\section{Boucksom's divisorial Zariski decompositions} 
\label{section:Zariski}

Let $(X,\omega)$  be a compact K\"ahler manifold.  
A cohomology class $\alpha\in H^{1,1}(X,\mathbb{R})$ is called pseudoeffective if it is the class of some positive closed current.  
Thanks to \cite[Theorem 3.12]{Boucksom04}, 
every pseudoeffective class $\alpha$ admits a unique divisorial Zariski decomposition 
$\alpha=Z(\alpha) + N(\alpha)$. 
Here the positive part $Z(\alpha)$ is a modified nef class, 
and the negative part $N(\alpha)$ is an effective divisor whose irreducible components form an exceptional family. 

According to \cite{Boucksom04}, a class $\beta$ is  modified nef if for any $\varepsilon>0$, 
there is a closed current $T\ge -\varepsilon \omega$ in the class $\beta$, such that $\nu(T,x) =0$ for any prime divisor $D\subseteq X$ and for any general point  $x\in D$.  
The cone   $\mathscr{E} \subseteq H^{1,1}(X,\mathbb{R})$ of pseudoeffective classes is called the pseudoeffective cone, 
and the one of modified nef classes is called the modified nef cone. 
A family of prime divisors is said to be exceptional, if the cone $\mathscr{F}$ generated by them meets the  modified nef cone  only at \{0\}.  
After \cite[Proposition 3.9]{Boucksom04}, 
the cone $\mathscr{F}$ is  a face of $\mathscr{E}$ in this case.
%on the boundary of $\mathscr{E}$.  
%In addition,  $\mathscr{E}$  is locally polyhedral around every point of $\mathscr{F}$ by \cite[Theorem 3.19]{Boucksom04}.   

%We will need the following version of negativity lemma. 

\begin{lemma}
\label{lemma:negativity-lemma} 
Let $f\colon X\to Y$ be a projective surjective morphism between compact K\"ahler manifolds. 
Let $E_1,...,E_r$ be a collection of prime  $f$-exceptional divisors.  
We assume that their sum is snc.  
Then the following properties hold.  
\begin{enumerate}
    \item For any non zero combination $D= a_1E_1+\cdots + a_rE_r$ with $a_1,...,a_r\ge 0$, 
    there is some $E_k$ such that $D|_{E_k}$ is not pseudoeffective. 
    \item $E_1,...,E_r$  form an exceptional family. 
    % In particular, the cone  $\mathscr{F}$  generated by  them   is on the boundary of  the pseudoeffective cone $\mathscr{E}$,   and  $\mathscr{E}$ is locally polyhedral around each point of $\mathscr{F}$.  
\end{enumerate}
\end{lemma}

We recall a that prime divisor $E$ in $X$ is called $f$-exceptional if $f(E)$ has codimension at least 2 in $Y$.   
Furthermore,  the precise meaning of $D|_{E_k}$ is the class 
$(\sum_{i=1}^k a_i c_1(E_i))|_{E_k}$. 

\begin{proof}
We first assume that the item (1) holds.  
Let $\mathscr{F}$ be the cone generated by the $E_1,...,E_r$. 
Then for any modified nef class $\gamma\in \mathscr{E}$, its restriction on any $E_k$ is pseudoeffective by \cite[Proposition 2.4]{Boucksom04}. 
Therefore, if  $\gamma  \in \mathscr{F}$, then it must be $0$. 
This proves the item (2).

It remains to prove the item (1).  
Let $Z\subseteq Y$ be an irreducible component of the image of the support of $D$. 
%the intersection of the centers of $E_1,...,E_r$. 
For any point $y\in Z$, we take an open Stein neighborhood $V$ of $y$, 
so that the preimage $W=f^{-1}(V)$ is embedded into  $V\times \mathbb{P}^N$, where $N>0$ is an integer. 
Let $V_1 \subseteq V$ be a submanifold,  cut out transversally by general hyperplanes, 
such that $V_1 \cap Z$ is a singleton, that the preimage $W_1:= f^{-1}(V_1)$ is smooth, and that $D|_{W_1}$ is a snc divisor on $W_1$.  
Afterwards, we take general hyperplanes of $V\times \mathbb{P}^N$, and cut $W_1$ into a smooth surface $W_2$. 
Let $V_2=f(W_2)$. 
Then the natural morphism $W_2\to V_2$ is a generically finite   projective morphism 
between complex analytic surfaces.  
The negativity lemma implies that, there is some $k\in \{1,...,r\}$, 
such that  $E_k\cap W_2$ is a $1$-cycle $C$ which has negative intersection with $D|_{W_2}$. 
This is equivalent to that $C\cdot (D|_{E_k})<0$ if we regard $C$ as a $1$-cycle in $E_k$. 
We notice that, the  class in  $H^{n-2,n-2}(E_k,\mathbb{R})$  of the integration over  $C$,  
is equal to the class of a product of  semipositive $(1,1)$-forms. 
This implies that $D|_{E_k}$ is not pseudoeffective, 
and completes the proof of the lemma.   
\end{proof}

As a corollary, we obtain the following assertion.

\begin{lemma}
\label{lemma:exceptional-family-bimero} 
Let $f\colon Y\to X$ be a projective bimeromorphic morphism between compact K\"ahler manifolds. 
Let  $\mathfrak{A}$ be a  family of prime divisors in $X$, 
let $\mathfrak{B} = \{f_*^{-1}D \ | \  D\in \mathfrak{A}\}$,   
and let $\mathfrak{C}$ be the  union of $\mathfrak{B}$ and the set of $f$-exceptional prime divisors. 
Then $\mathfrak{A}$ is an exceptional family if and only if so is $\mathfrak{C}$. 
\end{lemma}

\begin{proof}
Let $E_1,...,E_r$ be the prime $f$-exceptional divisors.   
Up to blowing up $Y$, we may assume that their sum is a snc divisor. 
Let  $ D_1,...,D_k $ be the prime divisors in $\mathfrak{B}$.  
We first assume that $\mathfrak{A}$ is an exceptional family. 
Let $\gamma$ be a modified nef class in $Y$ such that 
\[
\gamma = \sum_{i=1}^r a_i \cdot c_1(E_i)  +  \sum_{j=1}^k b_j  \cdot c_1(D_j) 
\]
for some real numbers $a_i,b_j \ge 0$.   
Then  $f_*\gamma =  \sum_{j=1}^k b_j \cdot c_1(f_*D_j) $ is a modified nef class on $X$. 
By the assumption on $\mathfrak{A}$, we deduce that  $b_1=\cdots b_k = 0$.  
Thus $\gamma = \sum_{i=1}^r a_i \cdot  c_1(E_i)$. 
By Lemma \ref{lemma:negativity-lemma}, we obtain that $\gamma=0$.

Suppose that $\mathfrak{C}$ is an exceptional family. 
Let $\eta $ be a modified nef class on $X$ such that 
\[
\eta =  \sum_{j=1}^k b_j \cdot c_1(f_*D_j) 
\]
for some real numbers $b_j \ge 0$. 
Then $f^*\eta$ is a pseudoeffective class and 
\[
f^*\eta = \sum_{i=1}^r a_i \cdot  c_1(E_i)  +  \sum_{j=1}^k b_j \cdot  c_1(D_j)
\]
for some real numbers $a_i\ge 0$.   
Since $\mathcal{C}$ is an exceptional family, the previous equation is the   divisorial Zariski decomposition of $f^*\eta$.   
Since $\eta$ is modified nef, the negative part of   $ f^*\eta$ is supported in the $f$-exceptional locus.   
It follows that $b_1=\cdots =b_k = 0$, and hence $\eta =0$. 
This completes the proof of the lemma. 
\end{proof}

The following statement plays an important role in our proofs. 

\begin{lemma}
\label{lemma:convexity} 
Let $V=  H^{1,1}(X,\mathbb{R})$  for some compact K\"ahler manifold $(X,\omega)$,  
let $\mathscr{E} \subseteq V$ be the pseudoeffective cone, 
and let $x_1,...,x_k\in \mathscr{E}$ be an exceptional family of prime divisor classes.   
Let $y\in V\setminus \mathscr{E}$ such that  
$y+a_1x_1+\cdots + a_kx_k  \not\in \mathscr{E}$ for any real numbers $a_1,...,a_k$. 
Let  $D\subseteq V$ be a bounded set. 
Then there is a    real number $\varepsilon_0>0$, 
such that if $z\in D$ and $0\le \varepsilon \le \varepsilon_0$, 
then $y + \varepsilon z +a_1x_1+\cdots + a_kx_k \not\in \mathscr{E}$  
for any real numbers $a_1,...,a_k$.    
\end{lemma}

\begin{proof}
Let $\mathscr{F}$ be the cone generated by $x_1,...,x_k$. 
Then it is a face of $\mathscr{E}$ by \cite[Proposition 3.9]{Boucksom04}.  
Thus, if $H$ is  the linear subspace generated by $\mathscr{F}$, 
then we have $\mathscr{F}=H\cap \mathscr{E}$.  
Let $\mathscr{G}\subseteq \mathscr{E}$ be the subset of pseudoeffective classes  whose negative parts do not contain any of the $x_i$. 
Then $\mathscr{G}$ is a closed convex cone by \cite[Proposition 3.5]{Boucksom04},  
$\mathscr{E}$ is generated by $\mathscr{G}$ and $\mathscr{F}$,  
and $\mathscr{G}\cap H = \mathscr{G}\cap \mathscr{F} = \{0\}$.  

Let $f\colon V \to V/H$ be the natural morphism, 
let $\mathscr{G}'=f(\mathscr{G})$ 
and  $y'=f(y)$. 
%Then $\mathscr{E}'$ is a convex  cone in $V'$.   
We note that  the   slice  
\[
S=\{\alpha \in \mathscr{G} \ | \  \alpha\cdot \{ \omega^{\dim X -1}\} = 1 \}
\]
is compact convex, and  that $\mathscr{G} = \mathbb{R}_{\ge 0} \cdot S$. 
Hence $f(S)$ is compact convex and does not contain the origin. 
It follows that  $\mathscr{G}' = \mathbb{R}_{\ge 0} \cdot  f(S)$ and is a closed convex cone.  
The assumption implies that $(y+H ) \cap \mathscr{E} = \emptyset$.
Hence $y' \notin \mathscr{G}'$.  
Therefore, there is a linear hyperplane $L'\subseteq V'$ which separates $y'$ and $\mathscr{G}'$. % Let $L\subseteq V$ be the preimage of $L'$.  
%Then it is a hyperplane containing $\mathscr{F}$, and separates $y$ and $\mathscr{G}$. 
Hence there is a linear form $\varphi$ on $V$ such that  $\varphi(y)  < 0$,   $\varphi|_\mathscr{G} \ge 0$  and $\varphi|_H \equiv 0$.   
Since $\mathscr{E}$ is generated by $\mathscr{G}$ and $\mathscr{F}$, 
we deduce that $\varphi|_\mathscr{E} \ge 0$.  
Since $D$ is bounded, there is a real number $\varepsilon_0>0$, 
such that  $ \varphi(y+\varepsilon z)  < 0$  
for all $z\in D$ and all   $0\le \varepsilon \le \varepsilon_0$.   
Hence $(y+\varepsilon z+ H) \cap \mathscr{E} =\emptyset$.  
This completes the proof of the lemma. 
\end{proof}

\section{Non pseudoeffective torsion-free sheaves}

In this section, we introduce a notion of non pseudoeffective torsion-free coherent sheaves on a compact K\"ahler manifold.

\begin{defn}
\label{def-non-psef} 
Let $\mathcal{E}$ be a torsion-free coherent sheaf on a compact K\"ahler manifold $X$, let $\mathbb{P}(\mathcal{E}) = \Proj (\mathrm{Sym}^\bullet \mathcal{E})$, 
and let $\mathcal{O}_{\mathbb{P}(\mathcal{E})}(1)$ be the tautological line bundle on $\mathbb{P}(\mathcal{E})$.  
We denote by $P\subseteq \mathbb{P}(\mathcal{E})$ the unique irreducible component which dominates $X$.  
We say that $\mathcal{E}$ is non pseudoeffective if the following property holds.  
Let  $\rho \colon Y\to P$ be  any desingularization, 
and let $f \colon Y\to X$ be the natural morphism. 
For any class $\delta\in H^{1,1}(Y, \mathbb{R})$ which is a linear combination of $f$-exceptional divisors, 
the class $c_1(\rho^*\mathcal{O}_{\mathbb{P}(\mathcal{E})}(1))+ \delta$ is not pseudoeffective on $Y$.  
\end{defn}

If $\mathcal{E}$ has rank one, then it is non pseudoeffective if and only if $\mathcal{L}:=\mathcal{E}^{**}$ is not a pseudoeffective line bundle. 
To see this, it is sufficient to take some smooth proper bimeromorphic model $f\colon X'\to X$ so that $(f^*\mathcal{E}/\mathrm{torsion})$ is locally free.  
If $\mathcal{E}$ is locally free, then it is non pseudoeffective if and only if  $\mathcal{O}_{\mathbb{P}(\mathcal{E})}(1)$  is not pseudoeffective.  
If $X$ is projective and $\mathcal{E}$ is reflexive, 
then $\mathcal{E}$ is non pseudoeffective if and only if it is not pseudoeffective in the sense of \cite[Definition 2.1]{HoeringPeternell2019}.

\begin{defn}
\label{def-movable} 
Let $X$ be a compact K\"ahler manifold of dimension $n\ge 2$. 
%Let $\mathscr{E} \subseteq H^{1,1}(X,\mathbb{R})$ be the closed cone generated by pseudoeffective classes. 
A class $\alpha \in H^{n-1,n-1}(X,\mathbb{R})$ is called movable, if it has non negative intersections with pseudoeffective classes.   
Let $\mathcal{E}$ be a torsion-free sheaf on $X$. 
The slope $\mu_{\alpha}(\mathcal{E})$ of  $\mathcal{E}$ with respect to $\alpha$ is the number 
\[
\mu_{\alpha}(\mathcal{E}) = \frac{c_1(\mathcal{E}) \cdot \alpha}{\rank \mathcal{E}}. 
\]
We also define  
\[ \mu_{\alpha,max}(\mathcal{E}) = \sup \{ \mu_{\alpha}(\mathcal{G}) \  |  \  0\neq \mathcal{G} \subseteq \mathcal{E} \} \in \mathbb{R} \cup \{+\infty\},\] 
\[ \mu_{\alpha,min}(\mathcal{E}) = \inf \{ \mu_{\alpha}(\mathcal{Q}) \  |  \  0\neq \mathcal{Q} \mbox{ is a torsion-free quotient of } \mathcal{E} \} \in \mathbb{R} \cup \{-\infty\}.\] 
The torsion-free sheaf $\mathcal{E}$ is called $\alpha$-semistable 
if $\mu_{\alpha}(\mathcal{E}) = \mu_{\alpha,max}(\mathcal{E})$, 
and is called $\alpha$-stable if $\mu_{\alpha}(\mathcal{G})<\mu_{\alpha}(\mathcal{E})$ for any non zero coherent subsheaf $\mathcal{G}$ of smaller rank. 
\end{defn} 

The following lemma was communicated to us by an anonymous person. 
It implies that $\mu_{\alpha,max}(\mathcal{E})$ and $\mu_{\alpha,min}(\mathcal{E})$ are always real numbers.  
Furthermore,   Harder-Narasimhan filtrations with respect to $\alpha$   exist.

\begin{lemma}
\label{lemma:A-bounded-above-v2}
Let $  X $ be a  compact  K\"ahler manifold   of dimension $n\ge 2$, 
let $\mathcal{E}$ be a coherent torsion-free sheaf, 
and  let   $\alpha \in H^{n-1,n-1}(X,\mathbb{R})$  be a movable class.   
Then the following set of real numbers is bounded from above 
\[
\{    {c}_1(\mathcal{G}) \cdot \alpha \  | \  0\neq  \mathcal{G} \subseteq \mathcal{E}\}. 
\] 
\end{lemma}

\begin{proof}  
We first observe the following fact. 
If $\mathcal{F}'\subseteq \mathcal{F}$  are coherent torsion-free sheaves of the same rank, 
then $\mu_{\alpha}(\mathcal{F}') \le \mu_{\alpha}(\mathcal{F})$. 
Indeed, there is a natural injective, generically isomorphic, morphism $\det \mathcal{F}' \to \det \mathcal{F}$. 
Hence $c_1(\mathcal{F}) - c_1(\mathcal{F}')$ is a pseudoeffective class. 
The inequality then follows for $\alpha$ is a movable class.

For any  coherent sheaf $\mathcal{G}$, we set $m(\mathcal{G}):= {c}_1(\mathcal{G}) \cdot \alpha$. 
Let $\Sigma$ be the set of non zero  subsheaves of $\mathcal{E}$. 
We define a partial order $\le$ on $\Sigma$, such that $\mathcal{G}\le \mathcal{G}'$ if and only if 
\begin{enumerate}
    \item $\mathcal{G}\subseteq \mathcal{G}'$, 
    \item  and $m(\mathcal{G}) \le m(\mathcal{G}')$.  
\end{enumerate}
Let $\mathcal{F}$ be a maximal element of $(\Sigma, \le)$ which has minimal rank.  
Such an element always exists, since   $\Sigma$  satisfies the ascending chain condition with respect to the inclusion.
We will show that $m(\mathcal{F}) \ge  \sup_{\mathcal{G}\in \Sigma}m(\mathcal{G})$, which implies the lemma.

Assume by contradiction that there is some $\mathcal{G}\in \Sigma$ such that $m(\mathcal{G}) > m(\mathcal{F})$. 
If $\mathcal{G} \not\subseteq \mathcal{F}$, then $\mathcal{F}\subsetneq \mathcal{F} + \mathcal{G}$. 
It follows that $m(\mathcal{F} + \mathcal{G}) < m(\mathcal{F})$ for  $\mathcal{F}$ is maximal.  
There is an exact sequence of  coherent sheaves 
\[
0\to \mathcal{F} \cap \mathcal{G} \to \mathcal{F} \oplus \mathcal{G} \to \mathcal{F}+ \mathcal{G} \to 0, 
\]
and hence  $m( \mathcal{F} \cap \mathcal{G})+ m(\mathcal{F} + \mathcal{G}) = m( \mathcal{F}  )+ m( \mathcal{G})$. 
We deduce that 
\[m( \mathcal{F} \cap \mathcal{G}) >  m(\mathcal{G}) > m(\mathcal{F}). \]
Therefore, replacing $\mathcal{G}$ by $\mathcal{F} \cap \mathcal{G}$, 
we can assume that $\mathcal{G} \subseteq \mathcal{F}$.  
Furthermore, we may assume that $\mathcal{G}$ is maximal among all elements $\mathcal{H}$ of $\Sigma$, such that $\mathcal{H} \subseteq  \mathcal{F}$ and that  $m(\mathcal{H}) > m(\mathcal{F})$.  

By the first paragraph, we must have $\mathrm{rank}\, \mathcal{G} < \mathrm{rank}\, \mathcal{F}$. 
Let $\mathcal{G}'$ be a maximal element, among the elements of $\Sigma$ containing $\mathcal{G}$. 
Then  $\mathcal{G}'$ is also a maximal element of $\Sigma$, and  $m(\mathcal{G}')   \ge m(\mathcal{G}) > m(\mathcal{F})$. 
By the first paragraph, either 
   $\mathcal{G}'  \not\subseteq \mathcal{F}$ 
   or $\mathrm{rank}\, \mathcal{G}' < \mathrm{rank}\, \mathcal{F}$.  
Since $\mathcal{F}$ is a  maximal element of minimal rank, we deduce that $\mathcal{G}' \not\subseteq \mathcal{F}$. 
Thus $m(\mathcal{F}+\mathcal{G}') < m(\mathcal{F})$.    
By the same argument of the previous paragraph, we obtain that 
\[m(\mathcal{F}\cap \mathcal{G}') > m(\mathcal{G}') \ge m(\mathcal{G}). \]
This is a contradiction, 
for   $\mathcal{G} \subseteq \mathcal{F}\cap \mathcal{G}' \subseteq \mathcal{F}$, and $\mathcal{G}$ is chosen to be maximal at the end of  the previous paragraph. 
\end{proof}

The following results shows a relation between tensor products and maximal slopes.  
When $X$ is projective, 
the statement was proved in  \cite[Theorem 4.2]{GrebKebekusPeternell2016}.  
Our proof follows the same strategy, and we will only highlight the different part. 

\begin{lemma}
\label{lemma:tensor-max-slope} 
Let $X$ be a compact K\"ahler manifold of dimension $n$, 
and let $\alpha\in H^{n-1,n-1}(X,\mathbb{R})$ be a movable class. 
Assume that $\mathcal{E}$ and $\mathcal{F}$ are torsion-free sheaves on $X$. 
Then we have 
\[
\mu_{\alpha,max }((\mathcal{E}\otimes \mathcal{F})^{**}) = \mu_{\alpha, max }(\mathcal{E}) + \mu_{\alpha, \max}(\mathcal{F}).  
\]
\end{lemma}

\begin{proof}
By the same argument of \cite[Proposition 4.4]{GrebKebekusPeternell2016}, 
it is sufficient to show that if both $\mathcal{E}$ and $\mathcal{F}$ are  $\alpha$-stable, then the reflexive tensor product  $(\mathcal{E}\otimes \mathcal{F})^{**}$ is $\alpha$-semistable.      
We note that  \cite[Proposition 2.7]{GrebKebekusPeternell2016} still holds in our setting. 
Hence we may assume that both $\mathcal{E}$ and $\mathcal{F}$ are locally free by blowing up $X$. 
Let $\omega$ be a K\"ahler form on $X$ and let $\eta = \omega^{n-1} \in H^{n-1,n-1}(X,\mathbb{R})$. 
Then by the same argument of  \cite[Theorem 3.4]{GrebKebekusPeternell2016} (see also  \cite[Lemma 6.7]{GrebRossToma2016}),  
for all $\varepsilon>0$ small enough, 
both $\mathcal{E}$ and $\mathcal{F}$ are stable with respect to $\alpha+\varepsilon\cdot \eta$.  
We note that in our situation of   compact  K\"ahler manifolds, 
it is sufficient to apply \cite[Lemma 4.3]{Toma2016} in the place of Grothendieck's boundedness lemma for this argument.  
Hence we may assume that $\alpha$ is in the interior of the movable cone, 
and we can now apply the method of Toma  in \cite[Section 6]{CampanaPeternell2011} to conclude that $\mathcal{E}\otimes \mathcal{F}$ is $\alpha$-semistable.  
This completes the proof of the lemma. 
\end{proof}

%is purely analytic,  and can be applied to our situation directly.    We note that  \cite[Lemma 5.5]{CampanaPeternell2011} has an alternative proof, by using the method of \cite[Theorem 3.4]{GrebKebekusPeternell2016}  or of \cite[Lemma 6.7]{GrebRossToma2016}.     To  adapt this method in our setting of compact  K\"ahler manifolds,  it is sufficient to apply \cite[Lemma 4.3]{Toma2016} in the place of Grothendieck's boundedness lemma. 

The remainder of this section is devoted to the following proposition.  

\begin{prop}
\label{prop:slope-non-psef}
Let $\mathcal{E}$ be a torsion-free sheaf on a compact K\"ahler manifold $X$, let $\alpha\in H^{n-1,n-1}(X,\mathbb{R})$ be a movable class  in the sense of Definition \ref{def-movable}. 
Assume that $\mu_{\alpha, max} (\mathcal{E}) < 0$. 
Then $\mathcal{E}$ is non pseudoeffective, in the sense of Definition \ref{def-non-psef}. 
\end{prop}

A crucial step of the proof  is the following statement, which relates the positivity of $\mathcal{O}_{\mathbb{P}(\mathcal{E})}(1)$ and the one of $\mathcal{E}$.   
The main ingredient of its proof is the positivity of the direct images of adjoint relative canonical sheaves, 
developed in a series of works, \textit{e.g.} \cite{Siu2002} \cite{Berndtsson2009}, \cite{BerndtssonPaun2008}, \cite{BerndtssonPaun2010}, \cite{Bocki2013},  \cite{GuanZhou2015},  \cite{PaunTakayama2018}, 
\cite{HaconPopaSchnell2018}, etc.

\begin{lemma}
\label{lemma:positivity-direct-image} 
Let     $f\colon Y\to X$ a be projective  fibration of positive relative dimension between  compact K\"ahler manifolds. 
Let $\mathcal{L}'$ be a $f$-ample line bundle on $Y$,  
let $E$ be an effective $f$-vertical divisor, 
and let 
\[\mathcal{L}  = \mathcal{L}' \otimes \mathcal{O}_Y(E).   \]
Assume that $\mathcal{L}$ is pseudoeffective. 
Then for any K\"ahler form $\omega$ on $X$,  
there is a singular Hermitian metric $h$ on $\mathcal{L}$, 
and a sufficiently large  integer $m$, 
such that the direct image sheaf
\[
\mathcal{E}:= f_*(\omega_{Y/X} \otimes \mathcal{L}^{\otimes m} \otimes \mathcal{I}(h^m)) 
\]
%is non zero, and admits a singular Hermitian metric $h_{\mathcal{E}}$ whose curvature $\Theta$ satisfies  \[ \Theta +  \omega \otimes \mathrm{Id}_{\mathcal{E}} \ge 0. \]
is non zero. Moreover,  $c_1(\mathcal{E}) +  \rank \mathcal{E}    \cdot  \{m \omega\} $ is a pseudoeffective class.
\end{lemma}

\begin{proof}

Since $\mathcal{L}'$ is $f$-ample, 
there is a  smooth Hermitian metric $h'$ on $\mathcal{L}'$ with curvature form $\theta $, 
such that  $\varepsilon \cdot \theta + \frac{1}{2} f^*\omega$  is a K\"ahler form on $Y$   
for some $\varepsilon>0$ small enough.  
%Thus $\{ \varepsilon \cdot ( \theta  + 2[E]) +  \frac{1}{3} f^*\omega\}$ is a big class for $E$ is effective. 
%Since $\mathcal{L}$ is pseudoeffective, 
%we deduce that $\{(1-\varepsilon)\theta + [E]  + \frac{1}{3} f^*\omega \}$ is a big class. 
Let $h_0$ be a smooth Hermitian metric on $\mathcal{L}$ with curvature form $\eta$.   
Since $\mathcal{L}$ is pseudoeffective, and since $E$ is effective,   the class 
\[
 \{\eta -\varepsilon\cdot \theta  \} = (1-\varepsilon) \cdot c_1(\mathcal{L}) + \varepsilon \cdot \{[E]\} 
\] 
is pseudoeffective. 
Thanks to Demailly's regularization in Theorem \ref{thm:regularization}, 
there is a  quasi-psh function  $\varphi$ with analytic singularities such that 
\[
T := \eta - \varepsilon \cdot \theta + \frac{1}{2} f^*\omega   + \ddc \varphi    
\]
is a positive current  on $Y$. 
It follows that 
\[
\eta + \ddc \varphi = T + \varepsilon \cdot \theta - \frac{1}{2} f^*\omega
\ge \varepsilon \cdot \theta - \frac{1}{2} f^*\omega
\ge -f^*\omega, 
\]
where the last inequality follows for  $\varepsilon \cdot \theta + \frac{1}{2} f^*\omega$ is a K\"ahler form.  

Let $h =  e^{-2\varphi} \cdot h_0$ be a singular Hermitian metric on $\mathcal{L}$.  
Then its   curvature current $\Theta_h$ satisfies 
\begin{equation}\label{eqn:almost-positive-current}
\Theta_h =  \eta + \ddc \varphi \ge -f^*\omega. 
\end{equation}
Let $x\in X$ be a general point and let $Y_x$ be the fiber over $x$.  
Then $Y_x$ is smooth and  $\varphi|_{Y_x} \not\equiv -\infty$. 
Therefore,  the curvature current of $h|_{Y_x}$ is equal to the restriction $\Theta_h|_{Y_x}$.  
We note that $\varphi|_{Y_x}$ has analytic singularities as well, and we have 
\[
\Theta_h|_{Y_x} = T|_{Y_x}+ \varepsilon \cdot \theta|_{Y_x}   - \frac{1}{2}f^*\omega|_{Y_x}
= T|_{Y_x}+ \varepsilon \cdot \theta|_{Y_x} 
\ge \varepsilon \cdot  \theta|_{Y_x}.  
\]
We remark that $ \theta|_{Y_x}$   is a K\"ahler form on $Y_x$, whose cohomology class is equal to  $c_1(\mathcal{L}'|_{Y_x})=c_1(\mathcal{L}|_{Y_x})$.  
By Lemma \ref{lemma:non-vanishing} below, 
we deduce that 
\[H^0(Y_x, \omega_{Y/X}|_{Y_x} \otimes \mathcal{L}|_{Y_x}^{\otimes m} \otimes  \mathcal{I}(h|_{Y_x}^m) ) \neq \{0\}\] 
for some positive integer $m>0$.   
Since  $\mathcal{I}(h|_{Y_x}^m) = \mathcal{I}(h^m|_{Y_x}) \subseteq \mathcal{I}(h^m)\otimes \mathcal{O}_{Y_x}$ after the Ohsawa-Takegoshi extension theorem, 
we deduce that  the  direct image sheaf $\mathcal{E}$ in the statement is non zero.

Let $X'\subseteq X$ be the largest open subset where $\mathcal{E}$ is locally free. 
Then $X\setminus X'$ has codimension at least $2$ in $X$. 
Let $X^\circ \subseteq X$ be a dense Zariski open subset, 
on which the following three properties hold. 
The morphism $f$ is smooth; 
the direct image  
$f_*(\omega_{Y/X} \otimes \mathcal{L}^{\otimes m})$ is  locally free satisfying the base change property; 
and $\mathcal{E}$ is a subbundle of $f_*(\omega_{Y/X} \otimes \mathcal{L}^{\otimes m})$.  
We can endow the canonical $L^2$ metric $h_{\mathcal{E}}$ on $\mathcal{E}$ induced by $h^m$, 
as shown in  \cite[Section 22]{HaconPopaSchnell2018}.

%Zariski open $U'$ subset of  $U$ by \cite[Proposition 1.1]{Raufi2015}, see also  \cite[Proposition 25.1]{HaconPopaSchnell2018}. 

Locally, let $o\in X$ be an    arbitrary point,  
let $U\subseteq X$ be a neighborhood of $o$ and let $V=f^{-1}(U)$.   
We are free to shrink $U$ around $o$. 
In particular, we can assume that $\omega|_U=\ddc \psi$ for some smooth psh function on $U$.    
It follows from \eqref{eqn:almost-positive-current} that $\overline{h} := e^{-2\psi \circ f}\cdot h$ 
is  a positive singular Hermitian metric on $\mathcal{L}|_V$. 
Moreover, $\mathcal{I}(h^m)= \mathcal{I}(\overline{h}^m)$. 
By \cite[Theorem 21.1]{HaconPopaSchnell2018}, 
the canonical $L^2$ metric $\overline{h}_{\mathcal{E}}$ on $\mathcal{E}|_{U}$ induced by     $\overline{h}$   has semipositive curvature.  
We note that  $h_{\mathcal{E}} = e^{2m\psi} \cdot \overline{h}_{\mathcal{E}}  $ on $U$.  
Indeed, by definition, for a section $s$ of $\mathcal{E}$ on $U$, 
for almost every  point $x$ in $U\cap X^\circ$, we have 
\[
h_{\mathcal{E}}(s(x)) = \int_{Y_x} h(\sigma(x,s)) 
=  e^{2m\psi(x)} \int_{Y_x} h(\sigma(x,s)) \cdot e^{-2m\psi(x)}   
= e^{2m\psi(x)} \overline{h}_{\mathcal{E}}(s(x)), 
\]
where $\sigma(x,s)\in H^0(Y_x, \omega_{Y_x} \otimes  \mathcal{L}^{\otimes m}|_{Y_x})$ is the element induced by $s$.

By \cite[Proposition 1.1]{Raufi2015} (see also  \cite[Proposition 25.1]{HaconPopaSchnell2018}), 
$\overline{h}_{\mathcal{E}}$ determines a singular Hermitian metric $\det \overline{h}_{\mathcal{E}}$ on $\det \mathcal{E}|_{X' \cap U}$, 
with positive curvature current.   
Up to shrinking $U$, we can assume that  $\det \mathcal{E}$ is trivial on $U$,  
and  that the determinant metric $\det h_{\mathcal{E}}$ on $\det \mathcal{E}|_{
U\cap X'}$ is represented by $e^{-2\eta}$, 
where $\eta\colon U\cap X' \to [-\infty,+\infty]$ is a measurable function.  
Since  $e^{-2m\psi} \cdot h_{\mathcal{E}}   =  \overline{h}_{\mathcal{E}}$, 
by taking the determinant of this equation, 
we deduce that 
\begin{equation*} %\label{eqn:potential}
\eta + \rank \mathcal{E} \cdot  m\psi  
\end{equation*}
is a psh function  on  $U\cap X'$.  
Since $ U\setminus  X' $ has codimension at least $2$ in $U$, the previous function extends to a unique  psh function on $U$, by the second main theorem of \cite{GrauertRemmert1956}.   
Since $\psi$ is smooth, we deduce that $\eta$ extends to a unique quasi-psh function on $U$. 
In conclusion, the singular Hermitian metric $\det h_{\mathcal{E}}$ on $\det \mathcal{E}|_{X'}$ extends to a  singular Hermitian metric on $\det \mathcal{E}$ over $X$, 
and the sum of its  curvature current  with $\rank{\mathcal{E}} \cdot m \omega$ is a positive.   
This completes the proof of the lemma.  
\end{proof}

%In addition,   by  \cite[Proposition 23.3]{HaconPopaSchnell2018},  for any section $t$ of $\mathcal{E}^*$ on $U$,  $\log \overline{h}_{\mathcal{E}}^*(t) $ is bounded from above,  where $\overline{h}_{\mathcal{E}}^*$ is the dual metric of $\overline{h}_\mathcal{E}$.  
%Since $\mathcal{E}$ is coherent, we may assume that it is finitely generated by sections over $U$.  We first suppose that $\mathcal{E}$ is locally free at $o$.  Up to shrinking $U$, we may assume that $\mathcal{E}$ is free on $U$.  It follows that $\det \overline{h}_{\mathcal{E}}^*$ is  bounded from above.  Equivalently, the psh function of \eqref{eqn:potential} is bounded from above on $U$. Since $\psi$ is smooth,   by the first main theorem of \cite{GrauertRemmert1956},  %(see also \cite[Th\'eor\`eme 1.7]{Demailly1985}),  we deduce that $\eta$ extends to a  unique quasi-psh function on $U$. Now we assume that $\mathcal{E}$ is not locally free at $o$.  The previous argument shows that $\eta$ is a quasi-psh function on the largest open subset $U'\subseteq U$ where $\mathcal{E}$ is locally free.  

The following assertion was used in the previous proof. 

\begin{lemma}
\label{lemma:non-vanishing} 
Let $X$ be a projective manifold, 
let $\mathcal{L}$ be an ample  line bundle  on $X$,    
let $h_0$ be a smooth Hermitian metric on $\mathcal{L}$ 
such that its curvature form  $\omega$ is K\"ahler.  
Assume that  $h$ is a singular Hermitian  metric on $\mathcal{L}$, 
such that it has analytic singularities and that 
its curvature current satisfies $\Theta_h \ge \varepsilon \omega$ for some $\varepsilon >0$.  
Then for any invertible sheaf $\mathcal{M}$, there is a positive integer $m$ such that 
\[
H^0(X, \mathcal{M} \otimes \mathcal{L}^{\otimes m} \otimes \mathcal{I}(h^m)) \neq \{0\}. 
\]
\end{lemma}

\begin{proof}
There is a quasi-psh function $\varphi$ such that $h  = h_0e^{-2\varphi}$.  
By assumption, for a sufficiently divisible integer $d>0$, $h_0^{d-2} e^{-2d\varphi}$ is a positive singular Hermitian metric on $\mathcal{L}^{d-2}$. 
Since $\varphi$ has analytic singularities, the multiplier ideal sheaves $\mathcal{I}(k\varphi)$ for any number $k>0$ can be computed in a single birational modification $\rho \colon X'\to X$, 
see for example \cite[Remark 5.9]{Demaillybook2012}. 
In particular, $\varphi \circ \rho$ is locally of the shape $\alpha \cdot \log |g| + O(1)$,  for some rational function $g$ and some real number $\alpha \ge 0$. 
Moreover, the vanishing loci of $g$ form a snc divisor in $X'$, with irreducible components   $E_1,...,E_r$.   
There is a divisor  $D=\sum_{j=1}^r \lambda_j E_j$ for some real numbers $\lambda_j\ge 0$, such that locally $D = \alpha \cdot \div (g)$. 
Since  $h_0^{d-2} e^{-2d\varphi}$ is a positive singular Hermitian metric on $\mathcal{L}^{d-2}$, 
we obtain from Theorem \ref{thm-decomposition} that $(d-2)\rho^*L - dD$ is  a pseudoeffective divisor, 
where $L$ is an  ample divisor corresponding to $\mathcal{L}$.  
Hence $(d-1)\rho^*L- dD$ is a big divisor.  
Thus, if $H$ is a divisor corresponding to $\mathcal{M}$,  then for some integer $m$ sufficiently large, the divisor 
\[
md \rho^*L -  \lfloor md D \rfloor + \rho^*H 
=  \lceil m(d-1)\rho^*L- mdD \rceil  +  \rho^* (mL+H)
\]
is linearly equivalent to an effective divisor.  
This shows that 
\[
H^0( X,  \mathcal{M} \otimes \mathcal{L}^{md} \otimes \rho_*\mathcal{O}_{X'}( -\lfloor md D \rfloor)) \neq \{0\}. 
\]

Let $\rho_j$ be the discrepancy of $E_j$ over $X$ for $j=1,...,r$. 
Then the multiplier ideal sheaves $\mathcal{I}(md\varphi)$ is equal to  
$ \rho_* \mathcal{O}_{X'}( \sum_{j=1}^r \rho_j  - \lfloor  md D \rfloor ).$ 
Since $\rho_j\ge 0$, we deduce that $\rho_*\mathcal{O}_{X'}( -\lfloor md D \rfloor) \subseteq \mathcal{I}(h^{md})$.  
This completes the proof of the lemma.  
\end{proof}

Now we are ready to prove Proposition \ref{prop:slope-non-psef}.  

\begin{proof}[{Proof of Proposition \ref{prop:slope-non-psef}}]
Assume by contradiction that $\mathcal{E}$ is   pseudoeffective. 
Let $P\subseteq \mathbb{P}(\mathcal{E})$ be the unique irreducible component which dominates $X$,  let  $\rho \colon Y\to P$ be a desingularization, 
and let $f \colon Y\to X$ be the natural morphism.     
We may assume that $Y$ is obtain by blowing up $P$ at smooth centers contained in the singular locus. 
In particular, there is an effective $f$-exceptional  divisor  
$D$  and some integer $k>0$  such that $ \rho^*\mathcal{O}_{\mathbb{P}(\mathcal{E})}(k) 
 \otimes \mathcal{O}_Y(-D)$ is $f$-ample.    
By assumption, there is an effective $f$-exceptional divisor $E$ such that  
\[
 \mathcal{L} : = \rho^*\mathcal{O}_{\mathbb{P}(\mathcal{E})}(k) 
 \otimes \mathcal{O}_Y(E)
\]
is pseudoeffective.  
Then $\mathcal{L}$ satisfies the condition of Lemma  \ref{lemma:positivity-direct-image} with respect to the fibration $f$.

Let $r=\rank \mathcal{E}$.  
For any integer $m \ge r$, 
we have 
\[\left( f_*(\omega_{X/Y} \otimes \mathcal{L}^{\otimes m}) \right)^{**}
\cong \det \mathcal{E} \otimes (\mathrm{Sym}^{mk-r} \mathcal{E})^{**},    
\]
since both of them are reflexive and they are isomorphic over the   locus where $\mathcal{E}$ is locally free. 
Let $\omega$ be an arbitrary K\"ahler form on $X$. 
We apply Lemma \ref{lemma:positivity-direct-image} to the morphism $f$. 
There is some singular Hermitian metric $h$ on $\mathcal{L}$, 
there is some integer  $m\ge r$, 
such that the following direct image is   non zero,    
\[
\mathcal{F} :=  f_*(\omega_{Y/X} \otimes \mathcal{L}^{\otimes m} \otimes \mathcal{I}(h^m)).  
\]
In addition, we have $\mu_{\alpha}(\mathcal{F}) \ge - m \cdot \{\omega\} \cdot \alpha$. 
Since 
\[\mathcal{F} \subseteq 
f_*(\omega_{X/Y} \otimes \mathcal{L}^{\otimes m}) 
\subseteq ( \bigotimes^{mk} \mathcal{E})^{**}, \]
it follows that 
\[
\mu_{\alpha, max} ( ( \bigotimes^{mk}  \mathcal{E})^{**})  \ge 
\mu_{\alpha}(\mathcal{F}) \ge  -m \cdot  \{\omega\} \cdot \alpha. 
\]
By Lemma \ref{lemma:tensor-max-slope}, we have $\mu_{\alpha, max} (( \bigotimes^{mk}  \mathcal{E})^{**}) = mk\cdot  \mu_{\alpha, max} (\mathcal{E})$. 
We hence deduce that $\mu_{\alpha, max} (\mathcal{E}) \ge - k^{-1}\cdot \{\omega\} \cdot \alpha$. 
Since $\omega$ is arbitrary and $k$ is fixed, this contradicts that $\mu_{\alpha,max}(\mathcal{E})<0$,   
and the proposition is proved.  
\end{proof}

%We will need the following statement. 

%\begin{lemma}
%\label{lemma:sing-locus-formal-ext} 
%Let $X$ be a compact complex manifold, let $C\subseteq X$ be a compact irreducible submanifold,   and let  $S_0$ be a locally irreducible locally closed subvariety of $X$.  
%Assume that $S_0$ contains a dense Zariski open subset $C_0$ of $C$,  that $S_0$ is smooth around a Zariski open subset $C_1$ of $C_0$,  and that $S_0$ extends formally along $C$.  Then there is a closed analytic subset $Z$ of $X$,  such that $Z\cap C_0$ is equal to $(S_0)_{\sing} \cap C_0$.  
%\end{lemma}

%\begin{proof} 
%We may assume that $C_1 = C_0\setminus (S_0)_{\sing}$. 
%As shown in Proposition \ref{prop:equi-ext-blowup}, there is a coherent ideal sheaf $\mathcal{J}$  on $X$ which is generated by $\mathcal{I}_{S_0}$ and $\mathcal{I}_C^2$ around $C_1$. 
%\end{proof}

\section{Notion of formal extensions}  
\label{section:blowup}

In order to apply  Lemma \ref{lemma:blowup-lelong} for the study of restriction Lelong numbers, one needs to keep blowing up and taking strict transforms.    
Let $X$ be a compact complex analytic variety, let $C\subseteq X$ be a compact irreducible subvariety. 
%Assume that $C$ is not contained in the singular locus of $X$. 
Suppose  that  ${S_0}\subseteq X$ is an irreducible  locally closed subvariety  of smaller dimension,  
which contains a   Zariski open subset $C_0\neq \emptyset$ of $C$.    
Assume that there is a smooth Zariski open subset of $C_1$ of $C_0$ such that $X$ and $S_0$   are smooth around $C_1$. 
Let $\rho\colon X'\to X$ be the blowup of $X$ at  $C$,  
and let ${S'_0}$ be the strict transform of ${S_0}$ in $X'$.  
In other words, ${S'_0}$ is the unique irreducible component of $\rho^{-1}(S_0)$ which dominates $S_0$. 
Let $E$ be the $\rho$-exceptional locus.  
Then  the intersection   of ${S'_0}$ and $E$ is isomorphic to $C_1$ over $C_1$.  
Thus there is a unique irreducible component   $C'_0$ of ${S'_0} \cap E$ which dominates $C$. 
Let  $C'$ be its Zariski closure  in $X'$.    
To apply  Lemma \ref{lemma:blowup-lelong}, 
we need to continue  blowing up $C_0'$ inside $X'$. 
More precisely, we need a proper bimeromorphic morphism $\rho'\colon X''\to X'$ such that, 
over a neighborhood of a general point of  $C_0'$, it is  the blowup at $C_0'$. 
A   necessary condition for the existence of $\rho'$ is that 
$C_0'$ is an open subset of $C'$.   
For this reason, we introduce the following notion.   

%However, if $C$ is a compact submanifold of a  complex manifold $X$, and if ${S_0}$ is locally closed submanifold containing only a proper open subset of $C$,  then the intersection inside the blowup of $X$ at $C$,  of the strict transform of ${S_0}$ and the exceptional locus,  may not be algebraic. In this case, we can no longer blow up  over $C$ as in Lemma \ref{lemma:blowup-lelong}.  

\begin{defn}
\label{def:blowup} 
With the notation above, if  $C'$ is isomorphic to $C_1$ over $C_1$, 
then we say that  $({S_0},C)$ admits a blowup.  
Recursively, we say that $({S_0},C)$ admits $m$ blowups for an integer $m>1$,   
if $({S_0},C)$ admit a  blowup and if $({S'_0},C')$ admits $m-1$ blowups. 
We say that $({S_0},C)$ admits infinitely many blowups,  
if $({S_0},C)$ admits $m$ blowups for any positive integer $m$. 
\end{defn}

\begin{remark}
\label{rmk:blowup-ext}
We note that if  $(S_0,C)$ admits a blowup, then the conormal bundle $\mathcal{N}_{C_1/{S_0}}^*$, defined on $C_1$,  
extends along $C$ to a coherent quotient  sheaf $\mathcal{Q}$ of $\mathcal{I}_{C}/\mathcal{I}_C^2$,  where $\mathcal{I}_C$ is the ideal sheaf of $C$ in $X$.  
To see this, it is enough to pushforward the coherent sheaf $\mathcal{O}_{X'}(-E) \otimes \mathcal{O}_{C'}$  onto $X$, where $\mathcal{O}_{X'}(-E)$ is the canonical  invertible ideal sheaf of the $\rho$-exceptional locus.    
Conversely, if such a quotient $\mathcal{Q}$ exists, 
then $C'$ is isomorphic to the main component of $\mathbb{P}(\mathcal{Q})$, 
and hence $(S_0,C)$  admits a blowup.  
The existence of $\mathcal{Q}$ is equivalent to that the sheaf $\mathcal{O}_{S_0}/(\mathcal{I}_{C_1}|_{S_0})^2$,  defined around $C_1$,  extends to a coherent quotient sheaf  of $\mathcal{O}_X/\mathcal{I}_C^2$.  
Hence we introduce the following notion which appears to be more natural. 
We will show that it is equivalent to  Definition \ref{def:blowup}.  
However,  we only need the properties in Definition \ref{def:blowup}, 
but not the ones in Definition \ref{def:extension}, 
throughout this paper.  
\end{remark}

\begin{defn}
\label{def:extension}
With the notation above,  for any integer $m\ge 1$, 
we say that $S_0$ extends along $C$    with order $m$, 
if the coherent sheaf $\mathcal{O}_{S_0}/(\mathcal{I}_{C_1}|_{S_0})^{m+1}$ defined around  $C_1$  extends on $X$ as a coherent quotient sheaf  of $\mathcal{O}_X/\mathcal{I}_C^{m+1}$.  
Equivalently, if there is a coherent ideal sheaf  $\mathcal{J}$ on $X$, 
which is generated by $\mathcal{I}_{S_0}$ and $\mathcal{I}_{C_1}^{m+1}$ around $C_1$.   
If  this is the case, we  say that $\mathcal{J}$ induces an extension of $S_0$ along $C$  with order $m$.
If $S_0$   extends along $C$ with any order $m$, 
then we say that  ${S_0}$ extends formally  along $C$.   
\end{defn}

\begin{prop}
\label{prop:equi-ext-blowup} 
With the notation above, 
%Let $X$ be a compact complex analytic variety,  let $C$ be a closed irreducible  subvariety, and let $S_0$ be an irreducible locally closed subvariety.  Assume that  $S_0\cap C$ contains a Zariski open subset $C_0$ of $C$,  and that there is a Zariski open subset $C_1$ of $S_0\cap C$ such that $X$, $C$ and $S_0$ are all smooth around $C_0$.  
$S_0$ extends formally along $C$ if and only if $(S_0,C)$ admits infinitely many blowups.  
In particular, in this case, the normal bundle $\mathcal{N}_{C_1/S_0}$ defined   on $C_1$ extends to a coherent sheaf on $C$. 
\end{prop}

\begin{proof}  
We first assume that  $S_0$ extends formally along $C$. 
Then $(S_0,C)$ admits a blowup, as pointed out in Remark \ref{rmk:blowup-ext}. 
Let $X'\to X$ be the blowup at $C$ and let  $m>1$ be an integer. 
Then there is   a coherent ideal sheaf $\mathcal{J}$ on $X$, 
which is generated by $\mathcal{I}_{S_0}$ and   $\mathcal{O}_X/\mathcal{I}_C^{m+1}$ around $C_1$.  Let  $\mathcal{J}' = \mathcal{O}_{X'} \cdot  \rho^{-1}\mathcal{J}$. 
Since $m\ge 0$, we see that $\mathcal{J}'$ is contained in  $\mathcal{O}_{X'}(-E)$, 
where $\mathcal{O}_X'(-E)$ is the canonical invertible ideal sheaf of the $\rho$-exceptional locus. 
We can then define the coherent ideal sheaf $\mathcal{J}'':=\mathcal{J}'\cdot \mathcal{O}_{X'}(E)$, 
and the coherent ideal sheaf $\mathcal{K}'$ generated by $\mathcal{J}''$ and $\mathcal{I}_{C'}^{m}$. 
Then $\mathcal{K}'$ induces an extension of $S'_0$ along $C'$ with order $m-1$.  
Thus $S'_0$ extends formally along $C'$.  
By induction, we can prove that $(S_0,C)$ admits infinitely many blowups.

Conversely, we  assume that  $(S_0,C)$ admits infinitely many blowups.  
Let $m\ge 1$ be an integer and let $\rho\colon X'\to X$ be the composition of $m$ blow ups defined inductively as in Definition \ref{def:blowup}. 
Let $S'_0$ be the strict transform of $S$ in $X'$,  
let $E$ be the $\rho$-exceptional locus. 
There is a unique irreducible component    of  $S'_0\cap E$ which dominate $C$. 
Let $C'$ be its Zariski closure.  
Let $\mathcal{J'}= \mathcal{I}_{C'}^{m+1}$ and 
let $\mathcal{J}\subseteq \mathcal{O}_X$ be the pushforward of  $\mathcal{J}'$ on $X$.
Then it is a coherent ideal sheaf  after Grauert's theorem. 
We claim that $\mathcal{O}_X/\mathcal{J}$ extends $\mathcal{O}_{S_0}/(\mathcal{I}_{C_1}|_{S_0})^{m+1}$. 
To see this, it is enough to observe that if $g$ is a local holomorphic function on $X$ around a   point of $C_1$, 
then its pullback $\rho^*g$ vanishes along $C'$ with order at least $m+1$ if and only if $g|_{S_0}$ vanishes  along $C_1$ with order at least $m+1$ (see for example Lemma \ref{lemma:blowup-lelong} in the case when $\dim S_0 = \dim C+1$).   
This implies that $S_0$ extends formally along $C$. 
Finally, the last assertion follows from the discussion in Remark \ref{rmk:blowup-ext}. 
This completes the proof of the proposition. 
\end{proof}

We have the following observation.

\begin{lemma}
\label{lemma:admit-blowup-bimero}
Let $X$, $C$, $S_0$  be as above %in the beginning of the section,   
and let $f\colon X \dashrightarrow \overline{X}$ be a proper bimeromorphic map.   
In other words, $\overline{X}$ and $X$ share  a same proper bimeromorphic modification.  
Assume that   $f $ is an isomorphism around general points of  $C$.   
Up to shrinking $S_0$, we may assume that it is contained in the isomorphic locus of $f$. 
Let $\overline{S}_0$ and $\overline{C}$ be the strict transform of $S_0$ and $C$ in $\overline{X}$. 
Then  $(S_0,C)$ admits  $m$ blowups if and only if so does $(\overline{S}_0,\overline{C})$.  
In particular, $S_0$ extends formally along $C$ if and only if $\overline{S}_0$ extends formally along $\overline{C}$. 
\end{lemma}

\begin{proof}   
It it enough to prove the case when $m=1$,  as for larger integers $m$, 
it is sufficient to repeat the same argument.    Let $\overline{\rho}\colon \overline{X}'\to \overline{X}$ and $\rho\colon X'\to X$ be the blowups at $\overline{C}$ and $C$ respectively,  let $\overline{E}$ and $E$ be the exceptional divisors.   Let  $D$  (respectively $\overline{D}$) be the intersection of the strict transform of ${S_0}$ in $X'$ and $E$  (respectively of the strict transform of $\overline{S}_0$ in $\overline{X}'$ and $\overline{E}$).   Then there is a natural proper bimeromorphic map $X' \dashrightarrow \overline{X}'$,  which is isomorphic  from a neighborhood of a general point in  $D$ to a neighborhood of a general point in $\overline{D}$. By passing to a common proper bimeromorphic modification of $\overline{X}'$ and $X'$,    we see that the Zariski closure  of $\overline{D}$ is equal to the strict transform of the Zariski closure  of $D$.  This proves the lemma.  
\end{proof}

The next lemma implies that,  
if $\mathcal{F}$ is a foliation on a compact complex manifold $X$, 
$C\subseteq X\times X$ is the diagonal and ${S_0}$ is the analytic graph of $\mathcal{F}$, 
%then $(S_0,C)$ admits infinitely many blowups. 
then  ${S_0}$ extends formally along $C$ in $X\times X$.

\begin{lemma}
\label{lemma:foliation-admit-blowup} 
Let $X$, $C$, $S_0$, $C_0$ and $C_1$ be as in the beginning of the section.  %we assume that $C$ is smooth. 
%that $C_0$ is a Zariski open subset of  $C$ whose complement has codimension at least $2$. 
We suppose further that $C_0=C_1$ and that  there is a foliation $\mathcal{F}$ on $X$, 
such that $\mathcal{F}$ is regular around $C_0$, 
that $\mathcal{F}$ is transversal to $C_0$, 
and that ${S_0}$ is the union of the local leaves of $\mathcal{F}$ passing through the points of  $C_0$.  
Then ${S_0}$ extends formally along $C$.  
\end{lemma}

We note that ${S_0}$ is called a graphic neighborhood of $C_0$   with respect to the foliation $\mathcal{F}$ in \cite{BogomolovMcQuillan16}. 

\begin{proof}
We will prove by induction on $m$ that $({S_0},C)$ admits $m$ blowups. 
For $m=1$, we let $\rho\colon X'\to X$ be the blowup at $C$, 
let ${S'_0}$ be the strict transform of ${S_0}$ in $X'$, 
let $C_0'=\rho^{-1}(C_0) \cap {S'_0}$, 
let $C'$ be the Zariski closure of $C_0'$, 
and let $\mathcal{F}'$ be the foliation on $X'$ induced by $\mathcal{F}$.  
Then by Lemma \ref{lemma:singularity-foliation-blowup} below, $C_0'$ is contained in the singular locus $\mathcal{F}'_{\sing}$ of $\mathcal{F}'$.  
Hence so is $C'$. 
Furthermore,  $\mathcal{F}'_{\sing}$ is indeed equal to $C_0'$
over a neighborhood of $C_0$.   
This shows the case of $m=1$.

We assume by induction that $({S_0},C)$ admits $m-1$ blowups for some integer $m\ge 2$.  
Let $\rho\colon X'\to X$ be the composition of the blowups, 
let ${S'_0}$ be the strict transform of ${S_0}$ in $X'$, 
let $E'$ be the $\rho$-exceptional locus, 
let $C_0' = E'\cap {S'_0}$  
and let $C'$ be the Zariski closure of $C'$ in $X'$.    
We let $q\colon \widehat{X}\to X'$ be the blowup of $C'$ and let $p\colon \widehat{X} \to X $ be the natural morphism. 
Let $\widehat{S}_0$  be the strict transform of ${S_0}$ in $\widehat{X}$, 
and let $\widehat{C}$ be the Zariski closure of the intersection $\widehat{C}_0$ of $\widehat{S}_0$ 
and the $p$-exceptional locus.  
By  Lemma \ref{lemma:singularity-foliation-blowup}, $\widehat{C}$ is an irreducible component of the singular locus $\widehat{\mathcal{F}}_{\sing}$, 
where $\widehat{\mathcal{F}}$ is the foliation on $\widehat{X}$ induced by $\mathcal{F}$.  
Furthermore, Lemma \ref{lemma:singularity-foliation-blowup} also implies that ${S'_0}$ and $E'$ intersect transversally. 
Hence $\widehat{C}_0$ is a connected component of $\widehat{\mathcal{F}}_{\sing} \cap q^{-1}(C_0')$. 
It follows that the natural morphism $\widehat{C} \to C'$ is isomorphic over  $C_0'$.   
This implies that $({S_0},C)$ admits $m$ blowups, 
and completes the proof of the lemma.  
\end{proof}

\begin{lemma}
\label{lemma:singularity-foliation-blowup} 
Let $X=Y \times Z$ be a product of complex manifolds, 
let $\mathcal{F}$ be the foliation induced by the natural projection $f\colon X\to Z$.  
Assume that $C\subseteq X$ is an irreducible submanifold such that $f|_C$ is an isomorphism onto its image in $Z$. 
Let $S = f^{-1}(f(C))$. 
We  suppose that $\dim S < \dim X$. 
Let $p \colon X' \to X$ be a composition of blowups, 
such that the first blowup center is $C$ and each blowup center afterwards is the intersection of the strict transform of $S$ with the exceptional locus. 
Let $\mathcal{F}'$ be the foliation on $X'$ induced by $\mathcal{F}$. 
%, or equivalently induced by the natural morphism ${f}'\colon {X}' \to Z$.   
Then the singular locus  $\mathcal{F}'_{\sing}$ of $\mathcal{F}'$ is equal to the singular locus of the reduced analytic subspace $ {f'}^{-1}(f(C))$. 
 
In particular, the fibers of  $\mathcal{F}'_{\sing}$ over $C$  are isomorphic to a disjoint union of copies of projective bundles over projective spaces. 
In addition, the intersection of the strict transform of $S$ and the $p$-exceptional locus is transversal. 
It is a component of $\mathcal{F}'_{\sing}$, 
and  is isomorphic to $\mathbb{P}(\mathcal{N}_{C/S})$, 
where $\mathcal{N}_{C/S}$ is the normal bundle of $C$ inside $S$. 
\end{lemma}

\begin{proof} 
The problem is local around every point of $C$. 
Locally, the inclusions $C\subseteq S \subseteq X$ are isomorphic to $\mathbb{C}^e \subseteq \mathbb{C}^{m+e} \subseteq \mathbb{C}^{m+n+e}$, 
which is the product of $\{0\} \subseteq \mathbb{C}^{m}\subseteq \mathbb{C}^{m+n}$ by $\mathbb{C}^e$.  
Thus, it is sufficient to prove the case when $e=0$, that is, when $C$ is a point.  
%We assume that $f(C)$ is the origin of $Z  \subseteq \mathbb{C}^n$.   

We will prove  the following statement by induction on the number $k$ of blowups in $p$. 
Let $Z'\to Z$ be the blowup of $f(C)$.  % with exceptional divisor $D$.  
Then the meromorphic map $g'\colon X'\dashrightarrow Z'$ is an equidimensional morphism on $X'\setminus S'$, with normal crossing fibers .  
For any point $x\in C'$, there are appropriate  coordinates systems so that the meromorphic map $g'$ can be written in coordinates around $x$ as 
\begin{equation}
\label{eqn:expression-mero}
(a_1,...,a_n,t,u_2,...,u_m) \mapsto (t^k a_1,\frac{a_2}{a_1},..., \frac{a_n}{a_1}).
\end{equation} 
In addition, $S'=\{a_1=\cdots = a_n=0\}$, 
the $g'$-exceptional locus is defined by $t=0$,  
and $C'=\{a_1=\cdots = a_n = t=0\}$. 
The lemma follows from this statement, since the foliation $\mathcal{F'}$ on $X'$ is equal to the one induced by $g'$.

We first study the case when $k=1$. 
Let $(\alpha_1,...,\alpha_n, \beta_1,...,\beta_m)$ be coordinates of $X$ so that $C$ is the origin and $S = \{\alpha_1=\cdots \alpha_n=0\}$.  
Assume that  $x\in E\setminus S'$. 
Without loss of the generality, we may assume that $x$ is not in the strict transform of $\{\alpha_1=0\}$. 
Then, in a neighborhood of  $x$,  the morphism $p\colon X'\to X$ can be written in coordinates as 
\[
(t,a_2,...,a_n,u_1,...,u_m) \mapsto (t,ta_2,...,ta_n, tu_1,...,tu_m). 
\]
Similarly,  $Z'\to Z$ can be locally written   as $(t,a_2,...,a_n)\mapsto (t,ta_2,...,ta_n)$.  
Thus $g'$ can be written in coordinates around $x$ as 
\[
(t,a_2,...,a_n,u_1,...,u_n) \mapsto (t,a_2,...,a_n). 
\]
This proves the property on $g'|_{X'\setminus S'}$. 
Assume that  $x\in C'$. Then we may assume that $p$ is written as 
\[
(a_1,a_2,...,a_n,t,u_2,...,u_m) \mapsto (ta_1,ta_2,...,ta_n, t, tu_2,...,tu_m). 
\]
Thus $g'$ can be written as 
\[
(a_1,a_2,...,a_n,t,u_2,...,u_m) \mapsto (ta_1,\frac{a_2}{a_1},...,\frac{a_n}{a_1}). 
\]
This proves the case when $k=1$. 

Before we proceed further to the inductive step, 
we observe the following remark. 
With the expression of \eqref{eqn:expression-mero}, 
the coordinates  values of $g'$ is not defined if $a_1=0$. 
Nevertheless, there is a bimeromorphic change of coordinates on $Z'$, 
which can be written as 
\[
(t,a_2,...,a_n) \mapsto (\frac{1}{a_2},ta_2, \frac{a_3}{a_2},...,\frac{a_n}{a_2}). 
\]
If we compose it with  \eqref{eqn:expression-mero}, 
then $g'$ can be written in new coordinates as 
\[
(a_1,...,a_n,t,u_2,...,u_m) \mapsto (\frac{a_1}{a_2},ta_2, \frac{a_3}{a_2},...,\frac{a_n}{a_2}). 
\]

Assume that we have proved the statement for some integer $k\ge 1$. 
Let $g\colon X'\to X$ be the composition of $k$ blowups, 
let $q\colon X''\to X'$ be the blowup at $C'$ with exceptional divisor $E''$, 
let $S''$ be the strict transform of $S'$, 
and let $C''$ be the intersection of $S''$ with $E''$. 
Let  $x\in E''\setminus S''$.  
From the previous paragraph, we  can  assume without loss of the generality that  $q$ is written in coordinates around $x$ as 
\[
(t',a_2',...,a_n',u_1',...,u_m') \mapsto (t',t'a_2',...,t'a_n',t'u'_1,u'_2,...,u'_m). 
\]
It follows that $g'\circ q\colon X''\dashrightarrow Z'$ can be written in coordinates around $x$ as 
\[
(t',a_2',...,a_n',u_1',...,u_m') \mapsto (t'^{k+1} \cdot u'^k_1,a_2',...,a_n').
\]
This proves the property on $(g'\circ q)|_{X''\setminus S''}$. 
Let  $x\in C''$.  
We can  assume that  $q$ is written in coordinates around $x$ as 
\[
(a'_1,a_2',...,a_n',t',u_2'...,u_m') \mapsto (t'a_1,t'a_2',...,t'a_n',t',u'_2,...,u'_m)
\]
It follows that $g'\circ q $ can be written in coordinates around $x$ as 
\[
(a'_1,a_2',...,a_n',t',u_2'...,u_m') \mapsto (t'^{k+1} a'_1 ,\frac{a_2'}{a_1'},...,\frac{a_n'}{a_1'}).
\]
This completes the  induction, 
and the proof of the lemma.  
\end{proof}

\section{Quasi-psh functions with large restriction Lelong numbers} 

The objective of this section is to prove Theorem \ref{thm:Zariski-dense-germ}. 
We first prove some preparatory lemmas.  

%Throughout this section, we let $X$ be a compact K\"ahler manifold of dimension $n$, let $C\subseteq X$ be an irreducible closed submanifold of dimension  $d$, and let $S\subseteq X$ be a locally closed submanifold containing $C$.  

\begin{lemma} 
\label{lemma:upper-bound-lelong}
Let $(X,\omega)$ be a compact K\"ahler manifold,  let $Z \subseteq X$ be a Zariski closed subset,  
and let $x\in X\setminus Z$ be a point. 
Assume that the prime divisors contained in $Z$ form an exceptional family.       
Let  $D\subseteq H^{1,1}(X,\mathbb{R})$ be a bounded subset. 
Then there is a constant $K$ such that, 
for any class $\alpha \in D$,  
for any real divisor class $\delta$ with support in $Z$, 
and for any positive closed current $T$ in the class $\alpha+ \delta$, 
the Lelong number $\nu(T,x)$ is bounded from above by $K$.   
\end{lemma}

\begin{proof} 
Let $f\colon Y\to X $ be the blowup at $x$, and let $E\subseteq Y$ be the exceptional divisor.  
Let $\mathscr{F}\subseteq H^{1,1}(X,\mathbb{R})$ be the cone generated by effective divisors supported in $Z$.   
Assume  that the lemma does not hold.  
Then for any number $N>0$, 
there is a closed positive current  $T_1$ in the class $\alpha+ \delta$ with $\alpha \in D $ and $ \delta \in \mathscr{F}$,   
such that   $\nu(T_1,x)>2N$.  
By Demailly's regularization in Theorem \ref{thm:regularization}, there is a closed positive current $T$ with class in $\alpha + \delta +\{\omega\}$ such that  $\nu(T,x) >N$ and that $T$ has analytic singularities.  
Thus, locally around $x$, we have $T=\ddc \psi$ for some  psh function $\psi$. 
Moreover, 
\[
\psi  = \frac{c}{2} \log (|g_1|^2+\cdots + |g_k|^2) + v,  
\]
where $c>0$ is a real number,   
$g_1,...,g_k$ are holomorphic functions, 
and $v$ is a  bounded function. 
Then $\nu(T,x)$ is equal to $c$ multiplied by the minimal vanishing order $m$ at $x$  of the  functions $g_1,...,g_k$. 

The pullback $f^*T$ is a closed positive current on $Y$, with cohomology class $f^*(\alpha+\delta + \{\omega\})$.  
For a general point $y\in E$, the Lelong number $\nu(f^*T, y)$ is equal to $\nu(T,x)$. 
By Siu's decomposition in  Theorem \ref{thm-decomposition}, if $[E]$ is the $(1,1)$-current of integration over $E$, then we have 
\[
f^*T - N[E] \ge  f^*T-\nu(T,x)[E] \ge 0.
\]
Thus  the cohomology class $f^*(\alpha+\delta + \{ \omega\}) - N \cdot c_1(E)$ is pseudoeffective.

We note that, together with $E$,  
the prime  divisors  in $Y$ with centers in $Z$   form an exceptional family, see Lemma \ref{lemma:exceptional-family-bimero}.   
Let $\mathscr{F}_Y \subseteq H^{1,1}(Y,\mathbb{R})$ be the cone generated by the effective divisors with centers contained in $Z$.   
Then  $-c_1(E)+\delta'$ is not pseudoeffective for any class $\delta' \in \mathscr{F}_Y$. 
%Moreover, the pseudoeffective cone $\mathscr{E}_Y$ is locally polyhedral around  every point of $\mathscr{F}_Y$.  
Since the image of $D$ under  $f^*$ is again bounded in $H^{1,1}(Y, \mathbb{R})$, 
we deduce from Lemma \ref{lemma:convexity} that there is a number  $K_0 > 0$,  such that 
\[
-K \cdot c_1( E)  + \delta' + f^*(\alpha'+\{\omega\})
\]
is not pseudoeffective, for any $K > K_0$, any class $\delta' \in \mathscr{F}_Y$, and any class $\alpha'\in D$. 
Since $N$ can be arbitrarily large, 
and since $f^*\delta\in \mathscr{F}_Y$, 
we obtain a contradiction to  the previous paragraph. 
This completes the proof of the lemma.
\end{proof}

The following key lemma is an application of Demailly's mass concentration.

\begin{lemma}
\label{lemma:con-mass}
Let $(X,\omega)$ be a compact K\"ahler manifold of dimension $n\ge 2$, 
and let $x\in X$ be a point.  
Let $z_1,...,z_n$ be local coordinate functions in a neighborhood $U$ of $x$ such that $x$ is the common zero of them.
Then for any number $\lambda >0$, 
there is an integer $m>0$,   
a $\omega$-psh function $\varphi$,  
such that, locally around $x$,  we have
\[
\varphi \le \frac{\lambda}{2m} \log (|z_1|^2 + |z_2|^{2m} + \cdots + |z_n|^{2m}) + O(1).
\]
\end{lemma}

\begin{proof}
This lemma is essentially the results in \cite[Section 6]{Demailly1993}.  
Although  in \cite[Section 6]{Demailly1993}, $X$ is assumed to be projective, 
and $\omega$ is assumed to be in the class of an ample line bundle, 
the same argument works for any compact K\"ahler manifold  $(X,\omega)$.  
In order to apply this method, it is sufficient that the mass of the Dirac measure $(\ddc \psi)^n$ at $x$ is smaller than the volume of $(X, \omega)$, 
where $\psi =  \frac{\lambda}{2m} \log (|z_1|^2 + |z_2|^{2m} + \cdots + |z_n|^{2m})$. 
Equivalently, it is enough that 
\[
(\frac{\lambda}{m})^n \cdot m^{n-1} <  \int_X \omega^n. 
\]
We note that the LHS above is equal to ${\lambda}^n  m^{-1}$. 
Hence the inequality   holds for some $m$ sufficiently large. 
This completes the proof of the lemma. 
\end{proof}

Now we are in position to prove Theorem \ref{thm:Zariski-dense-germ}.

\begin{proof}[{Proof of Theorem \ref{thm:Zariski-dense-germ}}] 
%We first remark that the property (1) already implies that $\varphi|_S \not\equiv -\infty$. 
%Indeed,  $\varphi^{-1}(-\infty) $ is equal to the cosupport of $\mathcal{J}$.  
%Therefore, it cannot contain  $S$, for $S$ is Zariski dense in $X$.   
We assume by contradiction that the theorem does not hold. 
Then there is a number $K$, 
such that  for any  $\omega$-psh function $\varphi$ satisfying the property (1)  of the theorem,  
we must have  $\nu(\varphi|_{S_0}, y) \le K$ for some point $y\in C_1$.     
Since $\varphi|_{S_0}$ has analytic singularities, 
this implies that  $\nu(\varphi|_{S_0}, y) \le K$ for general points $y\in C_1$.

Let $c$ be the codimension of $C$ in $X$. 
Let $x\in C_1$ be a point, and let $z_1,...,z_n$ be local coordinate functions in a neighborhood $U$ of $x$ in $X$, 
such that $x$ is the common zero of them, 
that $C\cap U$ is defined by $z_1=\cdots =z_c =0$, 
and that ${S_0}\cap U$ is defined by $z_1 = \cdots = z_{c-1} =0$.   
Let $\mu >0$ be an arbitrary   number.   
Then by Lemma \ref{lemma:con-mass}, there is a positive integer $m_1$, 
and a $\omega$-psh function $\psi_1$,  
such that locally around $x$,  we have
\[
\psi_1 \le \frac{\mu}{2m_1} \log (|z_1|^2 + |z_2|^{2m_1} + \cdots + |z_n|^{2m_1}) + O(1).
\]
Let $\psi_2$ be a regularization of $\psi_1$ as in Theorem \ref{thm:regularization}, so that the following properties hold.  
The current $2\omega + \ddc \psi_2 $ is positive.  
There is some large enough integer $m_2$, such that the quasi-psh function $\psi_2$ can be locally written as 
\[
\frac{1}{2m_2} \log (|h_1|^2 +   \cdots + |h_r|^{2}) + O(1), 
\] 
where  $h_1,...,h_r$ are local generators of a coherent ideal sheaf $\mathcal{J}$ on $X$, 
and we have $\mathcal{J} \subseteq \mathcal{I}(m_2\psi_1)$.
Since ${S_0}$ is Zariski dense in $X$,  it is not contained in the cosupport of $\mathcal{J}$. 
Thus one of the holomorphic functions $h_i$ is not identically zero along ${S_0}$.

We set $\psi_0 = \frac{\mu}{2m_1} \log (|z_1|^2 + |z_2|^{2m_1} + \cdots |z_n|^{2m_1})$,   which is a psh function on $U$.  
Then, locally around $x$, we have $\mathcal{J} \subseteq \mathcal{I}(m_2\psi_1) \subseteq  \mathcal{I}(m_2\psi_0)$. 
We note that, for any number $e>0$, the multiplier ideal sheaf $\mathcal{I}(em_1\psi_0)$ is generated on $U$ by the monomials $z_1^{a_1}\cdots z_n^{a_n}$ such that 
\[
a_1+1 + \frac{a_2+1}{m_1} + \cdots + \frac{a_n+1}{m_1} >\mu \cdot e. 
\]
Since $z_1$ vanishes long ${S_0}$,  by letting $e=\frac{m_2}{m_1}$ and $a_1=0$ in the previous inequality, we obtain that for any holomorphic function $h$ in $\mathcal{I}(m_2\psi_0)(U)$, 
the vanishing order of $h|_{U\cap {S_0}}$ at $x$ is greater than 
\[
(  \mu \cdot \frac{m_2}{m_1}  -1) \cdot m_1 - (n-1)
  =  \mu m_2-m_1+ 1 -n.  
\]
Let $\psi_3=\frac{1}{2}\psi_2$.  
It follows that the Lelong number $\nu(\psi_3|_{S_0}, x)$ is greater than 
\[
\frac{1}{2} \cdot \frac{1}{m_2}   \cdot (\mu m_2-m_1+ 1 -n) = 
\frac{1}{2}\mu - \frac{m_1+ n - 1}{2m_2}. 
\]
By Theorem \ref{thm:regularization},  $m_2$ can be chosen to be arbitrarily large once $\mu$ and $m_1$ are fixed.  
Thus we may assume that the previous number is at least $\frac{1}{3} \mu$.

We note that $\psi_3$ is a $\omega$-psh function which satisfies the property   (1) of the theorem.  
By assumption, for a general point $y\in C_1\cap U$, we have $\nu(\psi_3|_{S_0},y) \le K$. 
We write locally  $\psi_3 = \frac{1}{4m_2} \log (|h_1|^2 +   \cdots |h_r|^{2}) + O(1)$.  
Since ${S_0}$ extends formally along $C$,   
the normal bundle of $C_1$ in ${S_0}$ extends to a  line bundle  $\mathcal{N}$ on $C$, 
see Proposition \ref{prop:equi-ext-blowup}.  
Moreover, by Lemma \ref{lemma:blowup-lelong}, 
there is a composition of blowups $\rho \colon X' \to X$, 
such that the following properties hold. 
Over $C_1$, 
the strict transform ${S'_0}$ of ${S_0}$ in $X'$ meets the $\rho$-exceptional locus transversally  at a unique irreducible component  $E$. 
The minimal vanishing order $\zeta$ of $h_1\circ\rho$, ..., $h_r\circ \rho$ along $E$ satisfies 
\[
\frac{1}{2m_2} \cdot \zeta = \nu(\psi_3|_{S_0}, y) \le K,  
\]
where $y\in C_1$ is a general point.  
Let  $C'$ be  the Zariski closure of the unique component of  ${S'_0}\cap E$ which dominates $C$.   
Since $S_0$ is smooth around $C_1$, 
we deduce that $\rho(X'_{\sing})$ and  $\rho(C'_{\sing})$ are contained  in the set \[Z:=(C\setminus C_0)\cup (C_0\setminus C_1).  \] 
By blowing up $X'$ further at centers   contained in $Z$, 
we may assume that $X'$ and $C'$ are  smooth.  
Let  \[T':= \rho^*(\omega+\ddc \psi_3) - \frac{\zeta}{2m_2} [E],   \]
where $[E]$ is the current of integration over $E$. 
Then $T'$ is a closed positive current on $X'$ by  Theorem \ref{thm-decomposition}.   
In addition, a local potential of $T'$ is bounded   around general points of $C'$ by Lemma \ref{lemma:blowup-lelong}.

Let $C'_1$ is the preimage of $C_1$ in $C'$.  
Then $\rho|_{S'_0}\colon S'_0\to S_0 $ is an isomorphism around $C'_1$. 
Let $D$ be the union of the exceptional locus of $\rho|_{C'}$ and $(\rho|_{C'})^{-1}(C\setminus C_0)$.  
Then $C'\setminus D$ is an open subset of  $(\rho|_{C'})^{-1}(C_0)$, 
and it is isomorphic to its image in $C_0$. 
Hence $C'_1$ is a Zariski open subset of $C'\setminus D$,  
whose complement has codimension at least $2$. 
We note that 
\[\mathcal{O}_{X'}(E)|_{C'_1} \cong (\rho|_{C'_1})^*\mathcal{N}_{{C_1}/S_0}.\]
Since $S_0$  extends formally along $C$, there is a line bundle $\mathcal{N}$ on $C$ which extends  $\mathcal{N}_{{C_1}/S_0}$, see Proposition \ref{prop:equi-ext-blowup}.  
Hence by Lemma \ref{lemma:extend-class}, 
there is a divisor class $\delta'$ on $C'$ supported in $D$, 
such that 
\[c_1(E)|_{C'} = (\rho|_{C'})^*c_1(\mathcal{N}) + \delta'.\]
Let $T_C$ be the pushforward of $T'|_{C'}$ on $C$.   
Then it is closed positive and  its cohomology class satisfies  
\[ \{T_C\} = \{\omega|_C \} - \frac{1}{2m_2}\zeta \cdot  c_1(\mathcal{N}) + \delta,\] 
where  $\delta$ is a divisor class supported in $C\setminus C_0$.  
Let $x'\in C'_1$ be the point lying over $x$.  
Since ${S'_0}\to {S_0}$ is isomorphic over $x$, 
and since ${S'_0}$ and $E$ meet  transversally at $x'$,  
we can compute that 
\[
\nu(T_C,x) = 
\nu(T'|_{C'},x') \ge \nu(T'|_{{S'_0}}, x') 
= \nu(\psi_3|_{S_0}, x) - \frac{1}{2m_2}\zeta 
\ge \frac{\mu}{3} - K.  
\] 
Since $\mu$ is sufficiently large,  
and since $\frac{1}{2m_2}\zeta\in [0,K]$, 
we obtain a contradiction 
by applying Lemma \ref{lemma:upper-bound-lelong} on $C$. 
This completes the proof of the theorem. 
\end{proof}

\section{Foliations induced by meromorphic maps}  

We will  prove  Theorem \ref{thm:alg-germ-nonpsef-normal},   Theorem \ref{thm:foliation-psef-alg} 
and Theorem \ref{thm:BMQ}   in this subsection.

\subsection{Proof of Theorem \ref{thm:alg-germ-nonpsef-normal}}

This subsection focuses on Theorem \ref{thm:alg-germ-nonpsef-normal}.  
We first exclude the case when $S_0$ is Zariski dense in $X$.  
This also illustrates the idea of the complete proof of the theorem.

\begin{lemma}
\label{lemma:first-step-alg-germ} 
With the notation of Theorem \ref{thm:alg-germ-nonpsef-normal}, 
$S_0$ is not Zariski dense in $X$.  
\end{lemma}

\begin{proof} 
Up to shrinking ${S_0}$, we can assume that ${S_0}\cap C = C_0$. 
If we are in the situation of (4), 
then, by blowing up $X$ at $C$, and then blowing up further at centers  contained in $C\setminus C_0$, 
we can obtain a projective bimeromorphic morphism $\rho\colon X'\to X$ so that the following properties hold. 
Let $C_0'$ be the intersection  of the $\rho$-exceptional locus with the strict transform  ${S'_0}$ of ${S_0}$ in $X'$.      
If $C'$ is the Zariski closure of $C_0'$ in $X'$, 
then it is smooth and $C_0'$ is open in $C'$.   
We note that  ${S'_0}$ extends formally along $C'$   
by Lemma \ref{lemma:admit-blowup-bimero}. 
The   natural morphism $r\colon  C'\to C$  factors through $p\colon C' \to \mathbb{P}(\mathcal{N}^*)$. 
The conormal bundle $\mathcal{N}_{C_0'/{S'_0}}^*$ extends to the line bundle $\mathcal{N}'^* := p^*\mathcal{O}_{\mathbb{P}(\mathcal{N^*})}(1)$ on $C'$. 
The prime divisors in $C'$ contained in $C'\setminus C_0'$ are $r$-exceptional, and their sum is snc by the construction of $X'$.   
Hence by Lemma \ref{lemma:negativity-lemma}, 
they  form an exceptional family.  
By Definition \ref{def-non-psef}, $c_1(\mathcal{N}'^*)+\delta'$ is not pseudoeffective for any divisor class $\delta'$ supported in $C'\setminus C_0'$. 
In conclusion, up to  replacing $X$, ${S_0}$ and $C$ by $X'$, ${S'_0}$ and $C'$ respectively,  in the remainder of the proof,  
we can assume  in the situation of (3).

We assume by contradiction that ${S_0}$ is  Zariski dense   in $X$.   
For an arbitrary number $\lambda>0$, 
let $\varphi$ be a  $\omega$-psh function as in the conclusion of  Theorem \ref{thm:Zariski-dense-germ}.     
In particular,  $\varphi$ can be locally written as 
\[\varphi = \frac{1}{2m} \log  (|g_1|^2 +  \cdots |g_r|^{2}) + O(1)\]  
where  $g_1,...,g_r$ are local generators of  a coherent ideal sheaf $\mathcal{I}$ on $X$. 
By Lemma \ref{lemma:blowup-lelong}, 
there is a composition of blowups $\rho \colon X' \to X$, 
such that the following properties hold. 
The strict transform ${S'_0}$ of ${S_0}$ in $X'$ meets the $\rho$-exceptional divisor transversally  at a unique irreducible component  $E$. 
The minimal vanishing order $k$ of $g_1\circ\rho$, ..., $g_r\circ \rho$ along $E$ satisfies 
\[
\frac{k}{m} = \nu(\varphi|_{S_0}, y) > \lambda,  
\]
where $y\in C_0 $ is a general point. 
Let  $C'$ be  the Zariski closure of ${S'_0}\cap E$. 
Up to blowing up $X'$ further  at centers strictly contained in $C\setminus C_0$, 
we can assume that $C'$ and $X'$ are smooth.  
Let 
\[T':= \rho^*(\omega+\ddc \varphi) - \frac{k}{m} [E].  \]
Then $T'$ is a closed positive current by  Theorem \ref{thm-decomposition}.  
Moreover, a local potential of $T'$   is  bounded  around a  general point of $C'$ by the construction of Lemma \ref{lemma:blowup-lelong}.      
Let $T$ be the pushforward of $T'|_{C'}$ on $C$.   
Then it is closed and positive.  
There is a  divisor class $\delta'$  supported in $C\setminus C_0$, 
such that 
\[
\{T\} = \{\omega\} + \frac{k}{m}  \cdot  c_1(\mathcal{N}^*) + \delta' 
\]
We note that $\frac{k}{m} > \lambda$ can be arbitrarily large. 
By applying Lemma \ref{lemma:convexity} to the pseudoeffective cone of $C$, 
the condition (4) of the theorem implies that  $\{T\}$ is not pseudoeffective for a sufficiently large $\lambda$. 
This is  a contradiction  
\end{proof}

If the Zariski closure $M$ of $S_0$ in $X$ is not equal to $X$, 
then we intend to apply Theorem \ref{thm:Zariski-dense-germ} on $M$. 
An obstruction is that $M$ may be singular.  
Nevertheless,  the singular locus of $M$ does not contain ${S_0}$.  
It is a natural idea to take a desingularization $M'$ by blowing up $M$.  
We note that, if a center of the blowup intersects $S_0$ along $B$, 
and if $B$ is an effective divisor in $C_0$, 
then  $c_1(\mathcal{N}^*_{C_0/S_0})$ will increase with some positive multiple of $c_1(B)$ if we work  inside $M'$. 
As a consequence, we may lose the condition that $\mathcal{N}^*$ is non pseudoeffective. 
To deal with this problem, we need to take the exceptional locus of the desingularization into account, for the study of restriction  Lelong numbers.  
We will first prove some preparatory lemmas.

\begin{lemma}
\label{lemma:blowup-resolution} 
Let $X$ be a complex analytic variety, 
let $C\subseteq S\subseteq X$ be irreducible closed subvarieties. 
Assume that $\dim C = \dim S-1$, 
that $C$ and $S$ are smooth, and that $S$ is not contained in the singular locus of of $X$.   
Then there is a projective bimeromorphic morphism $\rho \colon X'\to X$ satisfying the following properties. 
The morphism $\rho$ is obtained by firstly blowing up $C$, and then blowing up the intersection of the strict transform of $S$ 
and the exceptional locus for several times. 
If $S'$ is the strict transform of $S$ in $X'$ and if $C'$ is the intersection of $S'$ and the $\rho$-exceptional locus, 
then the singular locus of $X'$ does not contain $C'$.   
\end{lemma}

\begin{proof} 
We may assume that $X$ is singular along $C$. The problem is local around general points in $C$.  Let $o$ be a general point of $C$,  so that the embedding dimension of $X$ at $o$ is minimal among  the points of $C$.  
By shrinking $X$ around $o$,  we  may embed $X$ into   $\mathbb{C}^{m+1+n}$ so that $o$ is mapped to the origin.  Moreover,  $m=\dim C$ and   $m+1+n$ is the embedding dimension of $X$ at $o$.  Since $S$ and $C$ are smooth,  we may assume there is a coordinates system 
$(x_1,...,x_m, t , a_1,...,a_n)$ of $\mathbb{C}^n$, 
such that $C$ is defined as $a_1=\cdots=a_n=t=0$ and $S$ is defined by $a_1 =\cdots =a_n = 0$.  
We write $\mathbf{a}=(a_1,...,a_n)$ and $\mathbf{x}=(x_1,...,x_m)$.  Then there are local holomorphic functions $ f_1,...,f_k$ in $(\mathbf{x},t,\mathbf{a})$, which define   $X$ inside $\mathbb{C}^{m+1+n}$.   
By assumption on the embedding dimension,  the Jacobi matrix 
%$\mathrm{D} \mathbf{f}$ of $\mathbf{f}$ 
of $f_1,...,f_k$   is identically zero along $C$.  

For any   holomorphic function $G\in \{f_1,...,f_k\}$,  up to shrinking $X$ around $o$, 
we can write  \[
G(\mathbf{x},t,\mathbf{a}) = G_1(\mathbf{x},t,\mathbf{a}) +  G_2(\mathbf{x},t,\mathbf{a}) + \cdots 
\]
where $G_i$ is a homogeneous polynomial in $\mathbf{a}$ of degree $i$, with coefficients as holomorphic functions in $(\mathbf{x},t)$.  
Since $X$ is smooth at a point $( \mathbf{x}, t, \mathbf{0})$ with $t\neq 0$, we deduce that there is some $G\in \{f_1,...,f_k\}$ such that $G_1$ is not identically zero. 
Since $G_1$ is identically zero along $C$, there is an integer $e\ge 1$  such that $t^e$ divides all the coefficients of $G_1$ and $t^{e+1}$ does not divide some coefficient of $G_1$.  
Let $\rho \colon  X'\to X$ be the composition of $e$ blowup as in the lemma.  It can be embedded into a composition of blowups $\eta \colon W\to \mathbb{C}^n$ of the same type.  
Particularly, $X'$ is the strict transform of $X$ in $W$.  Moreover, if $o'$ is the unique point in $C'$ lying over $o$,  then there is a  coordinates system  $(\mathbf{x}',t',\mathbf{a}')$ in a neighborhood of $o'$ in $W'$,  such that 
$S'=\{a'_1=\cdots =a'_n =0\}$, the $\rho$-exceptional locus is $\{t'=0\}$, 
and $\rho$ can be written as 
\[
(\mathbf{x}',t',\mathbf{a}') \mapsto ( \mathbf{x}', t', t'^e a'_1,...,t'^e a'_n). 
\]
Thus  the  function
\[ G' (\mathbf{x}',t',\mathbf{a}') := t'^{-2e} \cdot (G\circ \eta) (\mathbf{x}',t',\mathbf{a}') \] 
belongs to the ideal of $X'$ inside $W'$,  see for example \cite[Section 3.C]{Ou2024} for more explicit computations. 
We note that the linear part  $G'_1$ in $\mathbf{a}'$ of $G'$ is not divisible by $t'$. 
Hence, at a general point of $C'$ near $o'$,  the embedding dimension of $X'$ is strictly smaller than the one of $X$ at $o$.   Since  the embedding dimension is always positive,  by repeating this argument,  we obtain the proof of the lemma.
\end{proof}

We will use the  following two lemmas  to study the restriction Lelong numbers on the exceptional locus of a desingularization $M'\to M$.

\begin{lemma}
\label{lemma:negativity-exceptional-curve}
Let $o\in S$ be a germ of smooth surface and let $C$ be a smooth curve passing through $o$.  
Let $f\colon Y\to S$ be a projective bimeromorphic morphism which is isomorphic over the complement of $\{o\}$. 
Assume that $Y$ is smooth. 
Let $C_0$ be the strict transform of $C$ in $Y$, 
and let  $C_1,...,C_k$ be the irreducible components of the $f$-exceptional locus. 
%Assume that $C_1$ intersects the strict transform $C_0$ of $C$. 
Let $\theta$ be a smooth $(1,1)$-form on $Y$,  
let $(\varphi^j)$ be a sequence of $\theta$-psh functions with analytic singularities 
on $Y$,   
let $\nu^j$ be the Lelong number of $\varphi^j$ at general points of $C_0$, 
and let $\mu^j_i$ be the Lelong number of $\varphi^j$ at general points of $C_i$ for $i=1,...,k$. 
Assume that  $(\nu^j)$ is an increasing divergent sequence, 
and that $\lambda_i$ is the coefficient of $C_i$ in $f^*C$. 
Then for all $i=1,...,k$,  we have 
\[
\liminf_{j\to +\infty} \frac{\mu^j_i}{\nu^j} \ge \lambda_i. 
\]
\end{lemma}

\begin{proof}
Let $R^j= T^j - \sum_{i=1}^k \mu^j_i [C_i] - \nu^j[C_0]$.   
Then $R^j$ is a positive closed current  on $Y$  by Theorem \ref{thm-decomposition}, 
and its cohomology class is equal to 
\[\{\theta\}- \sum_{i=1}^k \mu^j_i c_1(C_i) - \nu^j c_1(C_0).\]   
Moreover, since $\varphi^j$ has analytic singularities,  
local potentials of $R^j$ are bounded around general points of $C_0,C_1,...,C_k$.   
Hence for all $i=1,...,k$,  the restriction on $C_i$ of the previous class  is positive.   
Since the intersection matrix $(C_i\cdot C_j)_{1\le i,j\le k}$ is definite negative, 
there is a class $\eta$ of the shape  $\eta=\sum_{i=1}^k b_i \cdot c_1(C_i)$, 
such that $\{\theta\}|_{C_i} = \eta|_{C_i}$ for all $i=1,...,k$. 
Furthermore, since $f^*c_1(C)$ is zero on $C_1,...,C_k$, 
we deduce that, for all $i=1,...,k$,  the restriction on $C_i$ of the class 
\[
 \frac{1}{\nu^j} \cdot  \eta- \sum_{i=1}^k (\frac{\mu^j_i}{\nu^j} - \lambda_i) \cdot c_1(C_i)    
\]
is positive.  
In other words, this class is $f$-nef.  
The negativity lemma implies that $\frac{\mu^j_i}{\nu^j} - \lambda_i \ge - \frac{b_i}{\nu^j}$. 
Since the sequence $(\nu^j)$ diverges, 
we can conclude  the   lemma. 
\end{proof}

\begin{lemma}
\label{lemma:reducible-C}  
Let $M$ be a complex manifold, 
let $\theta$ be a closed smooth real $(1,1)$-form,  
let $ U$ be an irreducible locally closed subvarieties of $M$, 
and let $C,D\subseteq U$ be closed irreducible   subvarieties. 
Assume that  $\dim U -1 = \dim C = \dim D$,  
that $C$  intersects $D$ transversally along a prime divisor $B$ of $C$, 
and that $U$ is smooth around   $B$.  

We suppose that there is a complex manifold  $V$  satisfying the following properties.  
\begin{enumerate}
%    \item  $U$  meets $B$ and is smooth around $ B\cap U$. 
%    \item  $C$ and $D\cap U$ intersect transversally along $B\cap U$ in $U$.
    \item  There is  a projective bimeromorphic morphism $f\colon U\to V$ such that
    % $V$ is smooth and  
    $f$ induces  an isomorphism from $C$ to a smooth divisor  $Z\subseteq V$.  
    \item  $f(D)=f(B)$ is a smooth divisor in $Z$, and $D$ is $\mathbb{P}^1$-bundle over  $f(D)$. 
\end{enumerate} 
Let $\lambda$ be the coefficient of $D$ in the Cartier divisor $f^*Z$ on $U$.  

Assume that we have a sequence $(\varphi^j)$ of $\theta$-psh functions with analytic singularities,  such that $\varphi^j$ is bounded around general points of $S$, 
that the generic Lelong numbers 
\[\nu^j:=\inf_{x\in C\cap S_0} \nu(\varphi^j|_{S_0}, x)\] form a divergent increasing sequence.   
If $\mu^j$ is the Lelong number of $\varphi^j|_{S_0}$ at general points of $  D$,  then  we have 
\[
\liminf_{j\to +\infty} \frac{\mu^j}{\nu^j} \ge \lambda.  
\] 
\end{lemma}
  
We note that $S_0$ is not assumed to be smooth in the lemma.  
Thus the definitions of $\nu^j$  and of $\mu^j$  only consider  the points  contained in the smooth locus of $S_0$.

\begin{proof}  
Let $\Gamma$ be a smooth local  surface in $V$, 
which intersects $f(D)$ transversally at a very general point,  
such that $\Gamma\cap Z$ is a smooth curve.  
Let $\Gamma' = f^{-1}(\Gamma)$. 
We may assume that, for all $j$, 
$\varphi^j$ is bounded at general points of $\Gamma'$, 
the Lelong number of $\varphi^j|_{\Gamma'}$ at general points of $\Gamma'\cap C$ is $\nu^j$, and the one at general points of $\Gamma'\cap D$ is $\mu^j$.   
Moreover, $\Gamma'\cap C$ and $\Gamma'\cap D$ form a snc divisor around the point $\Gamma'\cap B$ in $\Gamma'$.  
Let $h\colon Y\to \Gamma'$ be the minimal resolution.  
Then it is an isomorphism over the point $\Gamma'\cap B$.   
Let $\psi^j$ be the pullback of $\varphi^j$ on $Y$, 
and let $g\colon Y\to \Gamma$ be the natural morphism.  
Then $\mu^j$ and $\nu^j$ are equal to the generic Lelong numbers of $\psi^j$ along the strict transforms $h^{-1}_*(\Gamma'\cap D)$ and $h^{-1}_*(\Gamma'\cap C)$ respectively. 
Moreover, $\lambda$ is equal to the coefficient of $h^{-1}_*(\Gamma'\cap D)$ in the Cartier divisor $g^*(\Gamma\cap Z)$.
Then we can deduce the lemma by  applying Lemma \ref{lemma:negativity-exceptional-curve} to  the sequence $(\psi^j)$ on $Y$. 
\end{proof}

We can now finish the proof of  Theorem \ref{thm:alg-germ-nonpsef-normal}.

\begin{proof}[{Proof of Theorem \ref{thm:alg-germ-nonpsef-normal}}]  
As in Lemma \ref{lemma:first-step-alg-germ}, 
we may assume in the  situation of (3) and that $C_0=C\cap S_0$. 
Let $M$ be the Zariski closure of $S_0$ in $X$. 
Assume by contradiction that $\dim M>\dim S_0$. 
%By Lemma \ref{lemma:first-step-alg-germ}, we can assume that   $M\neq X$. 
We note that the singular locus $M_{\sing}$ of $M$ does not contain ${S_0}$, 
for $M_{\sing}$ is also a closed analytic subset of $X$. 
Let $\rho \colon  M'\to M$ be a projective bimeromorphic map obtained by successively applying the following constructions. 
\begin{enumerate}
\item[(i)] By blowing up $M$ at $C$ as shown in Lemma \ref{lemma:blowup-resolution}, 
we can assume that the singular locus $M_{\sing}$ of $M$ does not contain $C$.   
We note that, by  blowing up $M$ at centers contained in $C\setminus C_0$ afterwards,  
we still have the condition  that $C$ is smooth.   
The prime divisors in $C$ supported in $C\setminus C_0$ still form an exceptional family by Lemma \ref{lemma:exceptional-family-bimero}. 
We label by $B_1,...,B_d$ the irreducible components of $C \cap M_{\sing}$, which meet $C_0$  and which have dimensions equal to  $\dim C-1$. 
Let $W$ be an  open subset  of $S_0$, 
such that   it  intersects every element in $\{B_1,...,B_d\}$.  
We note that $W$ may not be connected and we are free to shrink $W$ throughout the proof. 

\item[(ii)] By blowing up $M$ at  centers contained in its singular locus, 
we obtain a desingularization  $\rho \colon M' \to M$. 
Thanks to the step (i), 
we see that the strict transform $C'$ of $C$ in $M'$ is well defined.
Let $S'_0$ and $W'$ be the strict transform of $S_0$ and $W$ respectively in $M'$.  
%We set $C'_0=S'_0\cap C'$. 

%\item[(iii)] Let $F'\subseteq \rho^{-1}(C)$ be the closed analytic subset of $M'$,  defined by applying Lemma \ref{lemma:ext-and-blowup} at every blowups in the previous step,  with respect to $W$ and the strict transforms  of $W$.  Up to shrinking $W$, we may assume that every irreducible component of $F' $ has dimension equal to $\dim C$.  Moreover, except $C'$,  the center in $C$ of every irreducible component of $F'$ belongs to $\{B_1,...,B_d\}$.    Let $W'$ be the  strict transforms of $W$.   Then $W'\to W$ is proper bimeromorphic.  
%and is  isomorphic over $W\setminus M_{\sing}$.  Thus, Up to shrinking $W$,   we may assume that   $W'$ is isomorphic to $W$ outside  $F'$. 

\item[(iii)] 
Since $C$ is smooth and since $C'\to C$ is projective bimeromorphic,  
the locus where $(\rho|_{C'})^{-1}$ is not isomorphic has codimension at least $2$ in $C$. 
Thus we can define the strict transforms  $B'_1,...,B'_d$  of $B_1,...,B_d$ in $C'$ respectively. 
Up to blowing up $M'$ at $B'_1,...,B'_d$, 
we may assume that $W'$ is smooth around $W'\cap B'_i$ for $i=1,..,d$, 
see Lemma \ref{lemma:blowup-resolution}.  
By blowing up at $B'_i$ further, we may assume there are unique irreducible subvarieties $D'_1,...,D'_d$ in $\rho^{-1}(C)\cap S'_0$, whose dimensions are equal to $\dim C$, 
such that $D'_i\cap W'$  intersects $C'\cap W'$ along $B'_i \cap W'$. 
In addition, the intersections are transversal, 
and   $\rho$ induces a $\mathbb{P}^1$-bundle structure on $D'_i\cap W'$ over $B_i\cap W$.  
We can still guarantee that $M'$ is smooth by blowing up at centers  over $B_i\subseteq M_{\sing}$.  
Finally, by  blowing up $M'$ further at centers contained in the singular locus of $C'$, 
we may assume that $C'$ is smooth.   
\end{enumerate}

We note that $S'_0$ and $S_0$ extend formally along $C'$ and $C_0$ respectively after the construction,  
see Lemma \ref{lemma:admit-blowup-bimero}. 
We also remark that  $\rho$ is an isomorphism over 
the smooth locus of $M$.   
Since $S_0$ is smooth,  the singular locus of $S'_0$ is contained in the $\rho$-exceptional locus. 
Hence the center $A$ in $C_0$ of  $(S'_0)_{\sing} \cap C'$  is contained in $C_0\cap M_{\sing}$. 
From the  step  (iii), we deduce that  $A$ does not contain  any element in  $\{B_1\cap C_0,...,B_d\cap C_0\}$.     
Thus, by the definition of $B_i$  in the step (i), 
$A$ is an analytic subset of $C_0$ of codimension at least $2$.  
Let $C_1$ be the largest open subset of $(C_0\setminus A)$ over which $\rho|_{C'}$ is isomorphic.  
Then $C_1$ is a Zariski open subset of $C_0$, whose complement has codimension at least $2$. 
Let $C_1'$ be the preimage of $C_1$ in $C'$. 
Then by construction,  we have 
\[C'_1 \subseteq  C' \setminus \big(\mathrm{Ex}(\rho|_{C'})\cup (\rho|_{C'})^{-1}(C\setminus C_0)\big). \] 
Since the RHS in the inclusion above is isomorphic to a Zariski open subset of $C_0$, 
we deduce that the complement inside it of $C'_1$ is an analytic subset of codimension at least $2$. 
Moreover, $S'_0$ is smooth around $C'_1$ for $C_1\cap A = \emptyset$. 
%there is a smooth  open neighborhood  $S'_1$ of  $C'_1$ in $S'_0$.  
%The complement of $C'_1$ in  is an analytic subset of codimension at least $2$.  
We note that the prime divisors in $C'$ contained in $\mathrm{Ex}(\rho|_{C'})\cup (\rho|_{C'})^{-1}(C\setminus C_0)$ form an exceptional family by Lemma \ref{lemma:exceptional-family-bimero}. 
\\

Let $\omega'$ be a K\"ahler form on $M'$. 
We can construct a sequence of $\omega'$-psh functions $(\varphi^j)$ with analytic singularities by applying Theorem \ref{thm:Zariski-dense-germ} to $M'$ with respect to  $(S'_0,C')$, 
so that the generic Lelong numbers 
\[
\nu^j:= \inf_{x\in C'_1} \nu(\varphi^j|_{S'_0}, x)
\]
form an increasing divergent sequence. 
%Let $D_1,...,D_t$ be the other irreducible components of $F'$,  so that $D_1,...,D_d$ are the components which intersect $C'\cap W'$.   Moreover, if  $B'_i=D_i\cap C'$,  then up to renumbering,  the center of $B'_i$ in $C$ is equal to $B_i$.  
Let $\mu_i^j$ be the Lelong number of $\varphi^j|_{S'_0}$ at general points of $D'_i$ for $i=1,...,d$.       
By applying Lemma \ref{lemma:reducible-C} on $M'$ around each $B'_i$ with  $U=W'$ and $V=W$,  
we deduce that 
\begin{equation}
\label{eqn:limsup} 
\limsup_{j\to +\infty} (\lambda_i-\frac{\mu^j_i}{\nu^j}) \le 0, 
\end{equation}
where $\lambda_i$ is the coefficient of $D'_i\cap W'$ in the Cartier divisor $(\rho|_{W'})^*(C\cap W)$ on $W'$.

Now we argue as in  Lemma \ref{lemma:first-step-alg-germ}. 
We fix some integer $j>0$, and    apply Lemma \ref{lemma:blowup-lelong} along $C'$.
%and then along each $D_i$ successively,  together with Lemma \ref{lemma:ext-and-blowup} with respect to $W$ and the strict transforms of $W$ at every blowup.  
We can hence obtain a projective bimeromorphic morphism $\tau \colon M'' \to M'$ so that the following properties hold.  
Up to blowing up $M''$ over its singular locus, 
we may assume that it is smooth.  
Let $S''_0$ be the strict transform of $S'_0$.   
There is a unique prime $\tau$-exceptional divisors $E$ in $M''$, 
such that $S''_0$ meets $E$ transversally over $C'_1$. 
Their intersection   is  open  in its Zariski closure $C''$ in $M''$. 
%such that the preimage of $C'_1$ in $S''_0$ is contained in $E$. 
%, 
%and the preimages of $D_1\cap S'_1,..., D_d\cap S'_1$ are contained in $E_1,...,E_d$ respectively. 
%These preimages are all openin their Zariski closures in $M''$, which we denote by $C'',D''_1,...,D''_d$.    
By blowing up further, we can assume that $C''$ is smooth.  
We remark that $C''\to C'$ is an isomorphism over $C'_1$, 
for $S'_0$ is smooth around  $C'_1$. 
The following current  
\[
R^j:= \tau^*(\omega' + \ddc \varphi^j) -  \nu^j [E]
\]
is closed and positive on $M''$.  
Moreover,  a local potential of $R^j$ is bounded around a general point of $C''$.  

We first restrict $R^j$ on $S''_0$. 
If $D''_i$ is the strict transform of  $D'_i$, 
then we have $\nu(R^j|_{S''_0}, x) = \nu(\varphi^j|_{S'_0}, \tau(x))$ for a general point $x$ in $D''_i$, since the $\tau$-exceptional locus does not contain $D''_i$. 
Now we restrict $R^j$ on $C''$. Let $B''_i$ be the strict transform  of $B'_i$ in $C''$.  
Then $B''_i\cap S''_0 \subseteq D''_i$. 
Hence for a general point $y\in B''_i\cap S''_0$, we have 
\[
\nu(R^j|_{C''}, y) \ge \nu(R^j|_{S''_0}, y) \ge \inf_{x\in D''_i} \nu(R^j|_{S''_0}, x) = \inf_{x\in D'_i} \nu(\varphi^j|_{S'_0}, x)  = \mu^j_i.  
\]
By applying  Theorem \ref{thm-decomposition} on $C''$, we deduce that 
\[
R^j|_{C''} - \sum_{i=1}^d\mu^j_i[B''_i]
\]
is a closed positive current on $C''$. 
Hence 
\begin{equation}
\label{eqn:psef-class-01}
\{\tau^*\omega'\}|_{C''} - \nu^j c_1(E)|_{C''} - \sum_{i=1}^d \mu^j_i c_1(B''_i)    
\end{equation}
is a pseudoeffective class on $C''$.  
We will push forward this class on $C'$.  
On the open subset $(\tau|_{C''})^{-1}(C'_1) \subseteq C''$,  
we have 
\[
\mathcal{O}_{M''}(-E)|_{C''} \cong  \mathcal{N}_{C''/S''_0}^*  \cong (\tau|_{C''})^* \mathcal{N}_{C'_1/S'_0}^*.  
\]
Let $\mathcal{N}'^*$ be a line bundle on $C'$ extending the conormal bundle $\mathcal{N}^*_{C'_1/S'_0}$ , see  Proposition \ref{prop:equi-ext-blowup}. 
By applying Lemma \ref{lemma:extend-class} and by  pushing forward \eqref{eqn:psef-class-01} on $C'$, 
we obtain that 
\begin{equation}
\label{eqn:psef-class-02}
\{\omega'|_{C'}\} + \nu^j c_1(\mathcal{N}'^*) - \sum_{i=1}^d \mu^j_i c_1(B_i') + {\delta^j}'
\end{equation}
is a pseudoeffective class on $C'$, 
where ${\delta^j}'$ is a divisor class supported in $C'\setminus C'_1$. 
In particular, the components of ${\delta^j}'$ are either exceptional over $C$ or with centers contained in $C\setminus C_0$. 

Next we will push forward to $C$.  
We first note that the conormal bundles satisfy 
\[\mathcal{N}^*_{C'_1/S'_0} \cong \mathcal{O}_{S'_0}(-C'_1)|_{C'_1} 
\mbox{ and }   \mathcal{N}^*_{C_1/S_0} \cong \mathcal{O}_{S'_0}(-C_1)|_{C_1}.\] 
We also recall that $\lambda_i$ is the coefficient of $D'_i\cap W'$ in the Cartier divisor $(\rho|_{W'})^*(C\cap W)$ on $W'$.  
Since $S'_0\to S_0$ is isomorphic over a neighborhood of $C_1\setminus M_{\sing}$, 
and since the irreducible components of $(C_1 \cap M_{\sing})$, which have dimensions equal to $\dim C -1$, 
are $ B_1\cap C_1,...,B_d\cap C_1$,  we deduce that  
\[
\mathcal{N}^*_{C'_1/S'_0} \cong  (\rho|_{S'_0})^*\mathcal{N}^*_{C_1/S_0}  \otimes \mathcal{O}_{C'}(\sum_{i=1}^d \lambda_i B'_i)
\] 
on $C'_1$.
Thus by Lemma \ref{lemma:extend-class}, we have 
\[
c_1(\mathcal{N}'^*) =   (\rho|_{C'})^* c_1(\mathcal{N}^* ) +  \sum_{i=1}^d \lambda_i \cdot  c_1(B_i') 
+ \gamma, 
\]
where  $\gamma$  is a divisor class supported in $C'\setminus C'_1$.  
%whose components are either exceptional over $C$ or with centers contained in $C\setminus C_0$, 
By pushing forward \eqref{eqn:psef-class-02} to $C$, and by letting $\theta=(\rho|_{C'})_*(\omega'|_{C'})$,  
we deduce that the following class  is pseudoeffective on $C$, 
\[
\frac{1}{\nu^j} \cdot  \{\theta\} +
c_1(\mathcal{N}^*) + \sum_{i=1}^d (\lambda_i -\frac{\mu^j_i}{\nu^j})_+ \cdot c_1(B_i) + \delta^j, 
\]
where $(\lambda_i -\frac{\mu^j_i}{\nu^j})_+ = \max\{0, \lambda_i -\frac{\mu^j_i}{\nu^j} \}$, 
and $\delta^j$ is a divisor class with support contained in $C\setminus C_0$.  
By applying Lemma \ref{lemma:convexity} on $C$, there is a real number $\varepsilon>0$, 
such that 
\[
c_1(\mathcal{N}^*) + \sum_{i=1}^d \alpha_i c_1(B_i) + \alpha_0 \{\theta\} +\delta
\]
is not pseudoeffective for any real numbers $0\le \alpha_0,...,\alpha_d\le \varepsilon$, 
and any  divisor class $\delta$ with support contained in $C\setminus C_0$. 
However, by \eqref{eqn:limsup}, for all $j$ sufficiently large, 
and for all $i=1,...,d$, 
we have 
\[
\frac{1}{\nu^j} \le \varepsilon \mbox{ and }  (\lambda_i -\frac{\mu^j_i}{\nu^j})_+ \le \varepsilon. 
\]
We obtain a contradiction, and complete the proof of the theorem. 
\end{proof}

\subsection{Proof of Theorem \ref{thm:foliation-psef-alg} 
and Theorem \ref{thm:BMQ}}
%We   need the following lemma for Theorem \ref{thm:foliation-psef-alg}.    

We first  recall the notion of  analytic graphs (see \cite[Section 2.1]{BogomolovMcQuillan16} and  \cite[Section 2.2.2]{Bost01}).  
Assume that $\mathcal{F}$ is a foliation on a compact complex manifold $X$. 
Let  $C\subseteq X\times X$ be  the diagonal and  let $C_0= C\cap (X_0\times X_0)$.  
There is a foliation $\mathcal{G} := p_1^{-1} 0 \cap p_2^{-1} \mathcal{F}$ on $X\times X$, 
where $p_1,p_2$ are the natural projections, and $0$ is the foliation by points on $X$.  
Then  the image under $p_2$ of a local leaf of $\mathcal{G}$ is a local leaf of $\mathcal{F}$. 
The analytic graph ${S_0}$ is the union of local leaves passing through the points of $C_0$.

\begin{lemma}
\label{lemma:foliation-graph} 
Let $X$ be a compact K\"ahler manifold and let $\mathcal{F}$ be a foliation on $X$.  
Let  $X_0\subseteq X$ be the regular locus of  $\mathcal{F}$.  
Then $\mathcal{F}$ is induced by a meromorphic map if its analytic graph ${S_0}\subseteq X_0\times X_0$ has the same dimension as its Zariski closure in $X\times X$. 
\end{lemma}

\begin{proof}
Let $M$ be the Zariski closure  of ${S_0}$ in $X\times X$, 
and let $g = p_1|_M \colon M\to X$ be the natural projection. 
Then for a general point $x\in X_0$, the subset $p_2(g^{-1}\{x\})$ in $X$ is the closure of the leaf of $\mathcal{F}$ passing through $x$. 
By \cite[Theorem 4.9.1]{BarletMagnusson2019},  
the morphism $g$ induces an analytic family of cycles $(C_t)_{t\in T}$ in $M$, 
parametrized by a dense Zariski open subset of $T$ of $X$, 
whose elements are general fibers of $g$. 
By \cite[Theorem 4.3.20]{BarletMagnusson2019}, the direct image 
$((p_2)_*C_t)_{t\in T}$ is an analytic family of cycles in  $X$. 
Hence there is an irreducible component $Z$ of the Barlet space of $X$,  
there is a Zariski closed subset $Y\subseteq Z$, 
such that some Zariski dense subset $Y^\circ$ of $Y$  parametrizes the cycles  in the family $((p_2)_*C_t)_{t\in T}$.   
Since $X$ is compact K\"ahler, $Z$ and hence $Y$ are compact by \cite[Corollary 4.2.76]{BarletMagnusson2019}.   
We note that a  cycle parametrized by a general point of $Y$ is an irreducible subvariety of $X$, 
and is tangent to the foliation $\mathcal{F}$. 
Since through any point of $X_0$, there is exactly one leaf of $\mathcal{F}$, 
we deduce that the underlying space of the universal cycle over $Y$ is bimeromorphic to $X$. 
It follows that there is a natural dominant meromorphic map $f$ from $X$ to $Y$, 
such that $\mathcal{F}$ is induced by $f$.  
\end{proof}

Now we can deduce Theorem \ref{thm:foliation-psef-alg} and Theorem \ref{thm:BMQ}.

\begin{proof}[{Proof of Theorem \ref{thm:foliation-psef-alg}}]  
We note that the item (2) follows from the item (1) and Proposition \ref{prop:slope-non-psef}.  
Thus it is sufficient to prove the item (1).
Let $X_0\subseteq X$ be the regular locus of  $\mathcal{F}$. 
Then $X\setminus X_0$ has  codimension at least 2.   
Let $\overline{X} = X\times X$, let  $ {C}\cong X$ be the diagonal of $X\times X$, let $ {C}_0= (X_0\times X_0) \cap    {C}$, 
and let $ {{S_0}}$ be the analytic graph  of $\mathcal{F}$.  
Then $ {{S_0}}$ is a locally closed submanifold of $X\times X$ containing $C_0$. 
Let $p_1,p_2$ be the natural projections from $\overline{X}$ to $X$, 
and let $ {\mathcal{G}} = p_2^{-1}\mathcal{F} \cap p_1^{-1} 0$, 
where $0$ is the foliation by points on $X$.  
Then ${\mathcal{G}}$ is foliation on $\overline{X}$. 
It is regular around ${C}_0$, 
and is transversal to  ${C}_0$.  
In addition,  ${{S_0}}$ is locally the union of leaves of passing through ${C}_0$.  
Thus we can apply  Lemma \ref{lemma:foliation-admit-blowup} to show that  
%$({{S_0}},{C})$ admits infinitely many blowups.  
${S_0}$ extends formally along $C$.  
The conormal bundle $\mathcal{N}_{{C}_0/{{S_0}}}^*$ is isomorphic to $\mathcal{F}^*|_{X_0}$. 
Since by assumption $\mathcal{F}^*$ is non pseudoeffective, 
we can deduce that ${{S_0}}$ has the same dimension as its Zariski closure in $\overline{X}$,  
by applying Theorem \ref{thm:alg-germ-nonpsef-normal} to the manifolds  
$ \overline{X}$, ${C}$ and ${{S_0}}$. 
By Lemma \ref{lemma:foliation-graph}, this completes the proof of the theorem. 
\end{proof}

\begin{proof}[Proof of Theorem \ref{thm:BMQ}] 
The proof is similar to the one of Theorem \ref{thm:foliation-psef-alg}. 
Up to blowing up $X$, we may assume that it is smooth. 
Let $Y = X\times W$, let $C \subseteq Y$ be the submanifold induced by the natural morphism $W\to X$, and let $\mathcal{F}_Y$ be the foliation $p_1^{-1}\mathcal{F}\cap p_2^{-1} 0$ on $Y$, where $p_1,p_2$ are the natural projections on $Y$ and $0$ is the foliation by points on $W$.     
Then $C$ is transversal to  $\mathcal{F}_Y$, 
and we can define  the graphic neighborhood ${S_0}$, 
which is the union of local leaves of $\mathcal{F}_Y$ passing through the points of $C$.  
By applying Theorem \ref{thm:alg-germ-nonpsef-normal} to $Y$, $C$ and ${S_0}$, 
we deduce that ${S_0}$ is has the same dimension as its Zariski closure $Z$ in $Y$. 
Let $g\colon Z\to W$ be the natural projection. 
Then $C$ induces  a section $\sigma\colon W\to Z$ of $g$. 
For any point $x\in V$, let $w\in W$ be a point lying over it. 
Then the irreducible component $D$ of $g^{-1}(\{w\})$ passing through $\sigma(w)$ has the same dimension as the rank of $\mathcal{F}$.  
In addition, $p_1(F)\subseteq X$ is the Zariski closure of the leaf of $\mathcal{F}$ passing through $x$. 
This completes the proof of the corollary. 
\end{proof}

\section{Uniruled compact K\"ahler manifolds}  

This section is devoted to Theorem \ref{thm:uniruledness}.  
The proof consists of several steps and we will divide the section into several subsections. 
We note that all complex manifolds we consider in this section are homotopy equivalent to finite CW-complexes. 
Hence their cohomology groups with coefficients in $\mathbb{Z}$ are finitely generated.

\subsection{Relative Bergman kernel metrics}

%\begin{lemma}
%\label{lemma:semicontinuity} 
%Let $f\colon X\to Y$ be a smooth morphism between complex manifolds,  and let $\mathcal{L}$ be a line bundle on $X$. 
%Assume that for any $y\in Y$, there is an integer $m_y>0$ such that $H^0(X_y, \mathcal{L}^{\otimes m_y}|_{X_y}) \neq \{0\}$.   
%Then there is an integer $m$ such that $H^0(X_y, \mathcal{L}^{\otimes m}|_{X_y}) \neq \{0\}$ for any $y\in Y$. 
%In addition, $f_*\mathcal{L}^{\otimes m} \neq 0$. 
%\end{lemma}

%\begin{proof}
%By upper semicontinuity, see \cite[Section 10.5.4]{GrauertRemmert84},  
%for any integer $k> 0$, the set 
%\[E_k:= \{ y\in   Y  \ | \   H^0(X_y, \mathcal{L}^{\otimes k}|_{X_y}) \neq \{0\}  \}\]
%is an analytic subset of $Y$. 
%The assumption implies that  $\bigcup_{k\ge 1} E_k =   Y$.  
%Thus we deduce from Baire's theorem that there is some $m>0$ such that $E_m=Y$.  
%In addition, there is a Zariski open subset $\emptyset \neq Y^\circ\subseteq Y$, 
%such that $f_*\mathcal{L}^{\otimes m}$ satisfies the basechange property, see \cite[Section 10.5.5]{GrauertRemmert84}.  
%It follows that $f_*\mathcal{L}^{\otimes m} \neq 0$. 
%\end{proof}

Let $f\colon X\to Y$ be a K\"ahler morphism between complex manifolds. 
The positivity of relative canonical line bundle $\omega_{X/Y}$ plays an important role in the  proof of Theorem \ref{thm:uniruledness}.  
For  morphisms between general complex manifolds, 
the theory has been developed  in a series of work,  \textit{e.g.} \cite{Siu2002},  \cite{Berndtsson2009}, \cite{BerndtssonPaun2008},  \cite{BerndtssonPaun2010},    \cite{Bocki2013},  \cite{GuanZhou2015},  \cite{Cao2017},  \cite{PaunTakayama2018},    
\cite{HaconPopaSchnell2018}, etc.   
We will need the following version of the positivity theorem. 

%We note that this lemma is well-known if $X$ and $Y$ are projective varieties. 
%The proof relies on two main ingredient. 
%The first one is the positivity of direct images,  
%see for example \cite[Theorem 4.13]{Campana2004}.  

%The second one is  Kawamata's cyclic covering trick, which reduces to the case  when $f$ has reduced fibers in codimension 1 of $Y$. 
%Since we assume $Y$ is projective in the lemma, the same construction applies. 

\begin{lemma}
\label{lemma:Bergman-indep}
Let $f\colon X\to Y $ be a   projective morphism between complex manifolds, 
let $\mathcal{L}$ be a line bundle on $X$, 
and let $h$ be a  smooth positive Hermitian metric on $\mathcal{L}$. 
Assume that there is an integer $m>0$, such that  $H^0(X_y,\omega_{X_y}^{\otimes m} \otimes \mathcal{L}|_{X_y}) \neq \{0\}$ for some general point $y\in Y$.  
Let $H$ be the $m$-relative Bergman Kernel metric  on $\omega_{X/Y}^{\otimes m} \otimes  \mathcal{L}$. 
Let $h^{-1}$ be the dual metric on $\mathcal{L}^*$,  
and let  $g=H\cdot h^{-1}$ be the Hermitian metric on  $\omega_{X/Y}^{\otimes m}$ induced by the canonical isomorphic morphism  $\mathcal{L}\otimes \mathcal{L}^*\to \mathcal{O}_X$.   
Then the curvature current $\Theta_g$ of $g$ exists, and satisfies 
\[
\Theta_g \ge -\Theta_h, 
\]
where $\Theta_h$ is the curvature form of $h$. 

%Let $Y^\circ\subseteq Y$ be a dense Zariski open subset, 
%let $X^\circ = f^{-1}(Y^\circ)$,   
%and let $\mathcal{L}'$ be a $f$-ample line bundle on $X$. 
Assume that $\psi\colon \mathcal{L}' \cong  \mathcal{L}$ is an isomorphism of line bundles on  $X$,  
that $h'$ is a positive  smooth Hermitian metric on $\mathcal{L}'$, and  that 
\[
h'= \psi^*h \cdot e^{-2\eta\circ f},  
\] 
where $\eta $ is a real smooth   function on $Y$. 
If  $g'$ is the Hermitian metric on  $\omega_{X/Y}^{\otimes m}$  defined in the same way as in the previous paragraph with respect to $(\mathcal{L}',h')$,  
then we have $g'=g$.  
\end{lemma}

\begin{proof}  
The first paragraph of the lemma is a consequence of \cite[Theorem 0.1]{BerndtssonPaun2010}
 (see also \cite[Theorem 3.5]{Cao2017}). 
For the second part of the lemma, 
%we may assume that $\psi$ is the identity automorphism of $\mathcal{L}$.  
%Let $\mathcal{E} = f_*(\omega_{X/Y}^{\otimes m} \otimes \mathcal{L})$. 
let  $Y^\circ\subseteq Y$ be a dense Zariski open subset, 
over which   $f$ is smooth,  
$f_* (\omega_{X/Y}\otimes \mathcal{L})|_{Y^\circ}$ is locally free and satisfies the base change property.  

The problem is  local on $Y$.  
Let $o\in Y$ be a point. 
We are free to shrink $Y$ around $o$. 
%In particular, we may assume that $\mathcal{E}$ is finitely generated by global sections,   
For any point $y\in Y^\circ$, %for local section $w$ of $\omega_{X/Y}^{\otimes m} \otimes  \mathcal{L}$, 
for any point $x\in X_y$, 
let $e$ be a local generator of $\omega_{X/Y}^{\otimes m}\otimes \mathcal{L}$ around $x$. 
In other words, $\omega_{X/Y}^{\otimes m}\otimes \mathcal{L} \cong \mathcal{O}_X \cdot e$ around $x$. 
Then $H(e(x))$ is defined as $B(e(x))^{-1}$, where 
\[
B(e(x)) = \sup \frac{|\frac{u(x)}{e(x)}|^2}{ (\int_{X_y} h(u)^{\frac{1}{m}} )^{{m}}}, 
\]
and the sup is taken over all non zero sections $u \in H^0(X_y,  \omega_{X_y}^{\otimes m} \otimes \mathcal{L}|_{X_y} )$,     
see for example  \cite[Section 3]{Cao2017}.  
If $e = e_1 \cdot e_2$, where $e_1$ is a local generator of  $\omega_{X/Y}^{\otimes m}$ and $e_2$ is a local generator of $\mathcal{L}$, 
then we have  
$g(e_1(x)) = B(e(x))^{-1} \cdot h(e_2(x))^{-1}. $ 

Let $e_2'$ be a local generator of $\mathcal{L}'$. 
Then $(\psi^*h)(e_2') = h(\psi(e_2'))$.
We note that $\psi(e_2') = \lambda \cdot  e_2$ for some local invertible holomorphic function $\lambda$ around $x$, as $\psi$ is isomorphic.   
Let $e'=e_1\cdot e_2'$ and $e'' = e_1\cdot \psi(e_2') = \lambda \cdot e$. 
Then $e'$ is a local generator of $\omega_{X/Y}^{\otimes m} \otimes \mathcal{L}'$, 
and $e'' $  is a local generator of $\omega_{X/Y}^{\otimes m} \otimes \mathcal{L}$.  
By abuse of notation, we also denote by $\psi$ the induced isomorphism from $\omega_{X/Y}^{\otimes m} \otimes \mathcal{L}'$ to $\omega_{X/Y}^{\otimes m} \otimes \mathcal{L}$. 
Then $\psi(e') = e''$. 
By applying the same calculus with $h'$, 
we obtain that 
\begin{eqnarray*}
B'(e'(x)) 
&=& \sup \frac{|\frac{u'(x)}{e'(x)}|^2}{ (\int_{X_y} h'(u')^{\frac{1}{m}} )^{{m}} }\\  
&=& \sup \frac{|\frac{\psi(u')(x)}{e''(x)}|^2}{ (\int_{X_y} e^{-\frac{2}{m}\eta(y)}  \cdot h(\psi ( u')) ^{\frac{1}{m}})^{{m}} }   \\
&=& e^{ 2\eta(y)} \cdot |\lambda(x)|^{-2} \cdot \sup  \frac{|\frac{\psi(u')(x)}{e(x)}|^2}{ (\int_{X_y}  h(\psi(u'))^{\frac{1}{m}}  )^{{m}}},  
\end{eqnarray*}
%\\&=& B(e(x)) \cdot e^{2\eta(y)}
where the sup is taken over all non zero sections $u' \in H^0(X_y,  \omega_{X_y}^{\otimes m} \otimes \mathcal{L}'|_{X_y})$.  
Since $\psi$ is an isomorphism, 
we obtain that 
\[
B'(e'(x))  = e^{2\eta(y)} \cdot |\lambda(x)|^{-2} \cdot B(e (x)). 
\]
It follows that 
\begin{eqnarray*}
g'(e_1(x)) &=& B'(e'(x))^{-1} \cdot h'(e_2'(x))^{-1} \\ 
&=& B(e(x))^{-1} \cdot  |\lambda(x)|^{2} \cdot e^{ -2\eta(y)} \cdot h(\psi(e_2')(x))^{-1} \cdot e^{2\eta(y)} \\ 
&=& B(e(x))^{-1} \cdot  |\lambda(x)|^{2} \cdot h( e_2(x))^{-1} \cdot  |\lambda(x)|^{-2}  \\ 
&=& g(e_1(x)).  
\end{eqnarray*}  
Since  $X\setminus f^{-1}(Y^\circ)$ has measure zero in $X$, 
and since the local potentials of $g$ and $g'$ are $\Theta_h$-psh functions, 
this implies the lemma. 
\end{proof}

\subsection{Locally projective morphisms induced by holomorphic 2-forms}

In the situation of Theorem \ref{thm:uniruledness}, if we assume that $H^{2,0}(X,\mathbb{C}) = \{0\}$, 
then $X$ is projective  by Kodaira's embedding theorem,  and we can conclude with \cite[Theorem 2.6]{BoucksomDemaillyPuaunPeternell2013}.   
Otherwise, we will show that every   non zero closed holomorphic $2$-form $\sigma$ on $X$ induces a meromorphic map $f\colon X\dashrightarrow Y$.

\begin{lemma}
\label{lemma:kernel-foliation}
Let $X$ be a compact K\"ahler manifold of dimension $n$ and let $\alpha\in H^{n-1,n-1}(X,\mathbb{R})$ be a movable class.  
Assume that $\mu_{\alpha}(T_X) >0$, and that there is some non zero closed holomorphic $2$-form $\sigma \in H^0(X,\Omega_X^2) \cong H^{2,0}(X,\mathbb{C})$.    
Let $\mathcal{K}$ be the kernel of  the  morphism 
\[ \lrcorner \sigma \colon T_X \to \Omega_X^1 \] 
defined by the contraction  with $\sigma$.  
Let $\mathcal{F} \subseteq \mathcal{K}$ be the maximal destabilizer with respect to $\alpha$.  
Then $\mathcal{F}$ is a  foliation induced by a meromorphic map $f\colon X\dashrightarrow Y$.   

Let $X'$ and $Y'$ be smooth K\"ahler bimeromorphic models of $X$ and $Y$ obtain by blowing up,  
so that the induced map $f'\colon X'\to Y'$ is a fibration. 
Then  the pullback  of $\{\sigma\}$ in $H^{2,0}(X',\mathbb{C})$ belongs to the image of  $H^{2,0}(Y',\mathbb{C})$. 
\end{lemma}

\begin{proof}  
Since $\sigma$ is closed, we see that $\mathcal{K}$ is a foliation by  Cartan calculus.  
We note that the generic rank of  $\lrcorner \sigma$ is an even number,  and denote it by $2a$.  
Then, by taking the   wedge product for $a$ times, $\sigma$ induces an injective   morphism  $  \mathcal{O}_X \to \Omega^{2a}_X$.  
Since $\omega_X$ is not pseudoeffective, we deduce that $2a<\dim X$, and hence $\mathcal{K}\neq 0$.  
Let  $\mathcal{Q} = T_X/\mathcal{K}$. 
Then the bilinear morphism $T_X \times T_X \to \mathcal{O}_X$ defined by $\sigma$ induces a bilinear morphism $\mathcal{Q} \times \mathcal{Q} \to \mathcal{O}_X$, 
which is generically non degenerate. 
Hence there is a generically isomorphic morphism $\mathcal{Q} \to \mathcal{Q}^*$, 
and it follows that $c_1(\mathcal{Q})\cdot \alpha \leq 0$ (we would like to thank St\'ephane Druel for pointing out this quick proof).

Thus  $ \mu_{\alpha}(\mathcal{F}) \ge \mu_{\alpha}(\mathcal{K}) \ge \mu_{\alpha}(T_X) > 0$.  
We will show that $\mathcal{F}$ is a foliation. 
We may assume that  $\mathcal{F}  \neq \mathcal{K}$.  
Since $2\mu_{\alpha,min}(\mathcal{F}) >  \mu_{\alpha,min}(\mathcal{F}) > \mu_{\alpha,max}(\mathcal{K}/\mathcal{F})$, 
we deduce that  any morphism from $\bigwedge^2\mathcal{F}$ to $\mathcal{K}/\mathcal{F}$ is zero.   
Since $\mathcal{K}$ is already a foliation, this implies that the morphism
$\bigwedge^2\mathcal{F} \to T_X/ \mathcal{F}$
induced by the Lie bracket is zero. 
Hence $\mathcal{F}$ is a foliation.  
By Theorem \ref{thm:foliation-psef-alg},  $\mathcal{F}$ is induced by a meromorphic map  $f\colon X\dashrightarrow Y$.

%We note that the generic rank of  $\lrcorner \sigma$ is an even number,  and denote it by $2a$.  
%Then, by taking the   wedge product for $a$ times, $\sigma$ induces an injective   morphism \[ \psi\colon \mathcal{O}_X \to \Omega^{2a}_X.  \]
%Since $\omega_X$ is not pseudoeffective, we deduce that $2a<\dim X$, and hence $\mathcal{K}\neq 0$.  
%Let  $\mathcal{Q} = T_X/\mathcal{K}$. Then  the following  natural morphism  induced by $\lrcorner \sigma$ 
%\[
%\bigwedge^{2a} T_X \to \Omega_X^{2a}
%\]
%factors through $ \bigwedge^{2a} \mathcal{Q}$.  

%We claim that its image is contained in the one of $\psi$.   
%Indeed, for any integer $k=1,...,a$, 
%we let  $u_1,...,u_{2k}$ be local holomorphic vector fields on $X$. 
%Then 
%\begin{eqnarray*}
%0&=&u_{2k}\lrcorner \cdots \lrcorner u_1 \lrcorner 
%(\sigma^{ a+k}) \\ 
%&=& (-1)^k \frac{(a+k)!}{a!} \cdot (u_1\lrcorner \sigma) \wedge \cdots \wedge %(u_{2k}\lrcorner \sigma) \wedge \sigma^{ a-k} \\ 
%&& + \sum_{m=1}^k \eta_m \wedge \sigma^{  a-k+m}, 
%\end{eqnarray*}
%where $\eta_m$ is a sum of holomorphic $(2k-2m)$-forms of the shape 
%\[
%\varphi \cdot  (v_1\lrcorner \sigma) \wedge \cdots \wedge (v_{2k-2m}\lrcorner \sigma), 
%\]
%$\varphi$ is a local holomorphic function, 
%and $v_1,...,v_{2k-2m} \in \{u_1,...,u_{2k}\}$.  
%By induction on $k$, we can show that the image of 
%$(\bigwedge^{2k} \mathcal{Q}) \wedge \sigma^{a-k}$ inside $\Omega_X^{2a}$ is contained in 
%the image of $\psi$. 
%This implies the  claim.  

We  note that  $Y$   belongs to the Fujiki's class, 
hence it  is bimeromorphic to a compact K\"ahler manifold, see 
\cite[Section IV.3]{Varouchas1989}.  
Therefore, we may assume that $Y'$ is a compact K\"ahler manifold. 
Let $\rho\colon X'\to X$ be the  natural morphism.   
By construction, the morphism  $T_{X'/Y'} \to \Omega_{X'}^1$ defined by the contraction with  $\rho^*\sigma$ is zero.    
Hence $\rho^*\sigma$ is the pullback of some closed holomorphic $2$-form on $Y$'. 
%It follows that $\rho^*\{\sigma\}$ is in the image of $H^{2,0}(Y',\mathbb{C})$. 
This completes the proof of the lemma.  
\end{proof}

In the next statement, 
we show that  $f'$ is locally projective if we choose the $2$-form $\sigma$ appropriately.

\begin{lemma}
\label{lemma:choose-2-form} 
With the notation of Lemma \ref{lemma:kernel-foliation}, we may choose $\sigma\in H^0(X,\Omega^2_X)$ so that the following property holds. 
There is a K\"ahler form $\omega'$ on $X'$
such that 
\[\{\omega' + \rho^* \sigma  + \rho^* \overline{\sigma} \} \in H^2(X',\mathbb{Z}), \]
where $\overline{\sigma}$ is the conjugate of $\sigma$.

In particular,  for any point $o'\in Y'$, 
there is an open neighborhood $U'\subseteq Y'$ of $o'$, 
a $f'$-ample line bundle $\mathcal{L}'$ on $V':=f'^{-1}(U')$,  
and a smooth Hermitian metric $h'$ on $\mathcal{L}'$, 
such that the curvature form of $h'$ satisfies 
$\Theta_{h'} =  \omega'|_{V'}$. 
\end{lemma}

\begin{proof}  
We choose $\sigma$ by  applying Lemma \ref{lemma:dense-Q} below. 
There is a closed real $2$-form $\Omega$  such that  its de Rham cohomology class
$\{\Omega\}  $ belongs to the image of  $  H^2(X, \mathbb{Q})$
and that $\Omega =  \omega + \sigma + \overline{\sigma}$, 
where $\omega$ is a K\"ahler form on $X$, and $\sigma$ is a closed $(2,0)$-form.   
Let $\rho\colon X'\to X$ be the natural morphism. 
We note that there is a smooth closed real $(1,1)$-form $\eta$ on $X'$, 
whose cohomology class is a linear  combination of $\rho$-exceptional divisors with rational coefficients, 
such that $ \rho^*\omega + \eta$ is a K\"ahler form.   

In the remainder of the proof, we replace $f\colon X \dashrightarrow Y$ by $f'\colon X'\to Y'$, $\omega$ by $ \rho^*\omega + \eta$,    $\Omega$ by $ \rho^* \Omega + \eta$, 
and $\sigma $ by $\rho^*\sigma$. 
Then $\{\Omega\}$ still belongs to $H^2(X,\mathbb{Q})$. 
By multiplying a positive integer, we assume further that $\{\Omega\} \in H^2(X,\mathbb{Z})$.    
By Lemma \ref{lemma:kernel-foliation}, there is a closed holomorphic $2$-form $\tau \in H^0(Y,\Omega^2_Y)$ such that $f^*\tau = \sigma$. 
Let  $\delta = \{ \tau + \overline{\tau} \} \in H^2(Y,\mathbb{R})$.
Then  
\[
\{ \sigma + \overline{\sigma}\} = f^*\delta \in H^2(X,\mathbb{R}). 
\]

Let $o\in Y$ be a point and let $U\subseteq Y$ be a  neighborhood of $o$ which is isomorphic to a polydisc.   
We set $V=f^{-1}(U)$. 
Then $\{\Omega|_V\} \in H^2(V, \mathbb{Z})$. 
We note that   $H^2(U,\mathbb{Z}) = \{0\}$.  
Thus we have   
\[ 
\{(\sigma+\overline{\sigma})|_V\} = (f^*\delta)|_V = f^*(\delta|_U) = 0 \in H^2(V,\mathbb{R}). 
\]  
It follows that 
\[\{\omega|_V\} = \{\Omega|_V\} \in H^2(V,\mathbb{Z}). \]  
Hence by \cite[Theorem 1.1]{Bingener1983},  there is a $f$-ample line bundle $\mathcal{L}$ on $V$ 
and a smooth Hermitian metric $h$ on $\mathcal{L}$,  
such that   the curvature form of $h$ satisfies 
$\Theta_{h} =  \omega|_{V}$.  
This completes the proof of the lemma. 
\end{proof}

\begin{remark}
\label{rmk:bimeromorphic} 
We remark that, in the previous lemma, 
once we have chosen $\sigma \in H^0(X,\Omega_X^2)$, 
the conclusion of the lemma still holds if we replace $X'$ and $Y'$ by higher blowups. 
It is sufficient to modify the K\"ahler form $\omega'$ by adding some smooth form whose cohomology class is a linear  combination of exceptional divisors with negative rational coefficients.  
\end{remark}

%We claim that the image of $\{\omega|_W\}= \{\Omega|_W\} $  in $H^2(W,\mathcal{O}_W)$  is zero.      Indeed, since $W$ is an open subset of $X$,  the following diagram is commutative,   
%\centerline{
%\xymatrix{
%H^2(X,\mathbb{Z})  \ar[d]  \ar[r]  & H^2(X,\mathcal{O}_X) \ar[d]   \\ 
%H^2(W,\mathbb{Z})    \ar[r]  & H^2(W,\mathcal{O}_W)   } }
%where the vertical arrows are restriction maps.
%Since the image of $\{\Omega\}$ in $H^2(X,\mathcal{O}_X)$ is the Dolbeault class $\{\tau\}_{\mathrm{Dol}}$ of $\tau$,  the image of $\{\omega|_W\}$  in $H^2(W,\mathcal{O}_W)$ is equal to $\{\tau\}_{\mathrm{Dol}}|_W$.  We note that $\{\tau\}_{\mathrm{Dol}} \in f^*(H^2(Y,\mathcal{O}_Y))$ by assumption.    Thus $\{\tau\}_{\mathrm{Dol}}|_W \in f^*(H^2(U,\mathcal{O}_U))$.  Since $U$ is a polydisc, we see that $H^2(U,\mathcal{O}_U) = 0$. This proves the claim. 

%Then for any point $y\in U$, we have $c_1(\mathcal{L}|_{W_y}) = \{ \omega|_{W_y}\} \in  H^{2}(W_y,\mathbb{R})$. 
%Since $\omega|_{W_y}$ is a K\"ahler form on the compact complex manifold $W_y$, 
%it follows that $\mathcal{L}|_{W_y}$ is ample. 
%Thus $\mathcal{L}$ is $f$-ample by \cite[Corollary 1.6]{Nakayama1987}.  

The following lemma was used in the previous argument. 

\begin{lemma}
\label{lemma:dense-Q} 
Let $W$ be a vector space of finite dimension  over $\mathbb{Q}$, and let $W_{\mathbb{R}}=W\otimes_{\mathbb{Q}} \mathbb{R}$.  
Assume that there is a decomposition $W_{\mathbb{R}}= W_1 \oplus W_2$ of real vector spaces.  
We equip   $W_{\mathbb{R}}$ with an Euclidean norm. 
Let $E\subseteq W_{1}$ be a  subset whose interior  is not empty.  
Then there is some element $a\in W$ such that $a=a_1+a_2$ with $a_1\in E$ and $a_2\in W_2$. 

In particular, we let   $W=H^2(X,\mathbb{Q})$ for some compact K\"ahler manifold, 
let $W_1 = H^{1,1}(X,\mathbb{R})$, 
$W_2 = (H^{2,0}(X,\mathbb{C}) \oplus H^{0,2}(X, \mathbb{C}) )\cap  H^{2}(X,\mathbb{R})$ 
and let $E\subseteq W_1$ be the K\"ahler cone.  
Then there is a closed real $2$-form $\Omega$  such that  its de Rham cohomology class
$\{\Omega\}  $ belongs to $  H^2(X, \mathbb{Q})$
and that $\Omega =  \omega + \sigma + \overline{\sigma}$, 
where $\omega$ is a K\"ahler form on $X$, and $\sigma$ is a closed $(2,0)$-form.  
\end{lemma}

\begin{proof}
Since $W_{\mathbb{R}}$ is finite dimensional, any two norms on it are equivalent. 
Hence  we may assume  that $W_1$ is orthogonal to $W_2$.   
By assumption, there is some $z_1 \in E$ and some $r>0$, 
such that the intersection with $W_1$ of the open ball $B(z_1,r)$ centered at $z_1$ with radius $r$ is contained in $E$. 
Since $W$ is dense in $W_{\mathbb{R}}$, there is some element $a \in B(z_1,r)\cap W$. 
If we decompose $a=a_1+ a_2$ with respect to the decomposition $W_{\mathbb{R}}= W_1 \oplus W_2$, then we have 
\[||a_1-z_1||^2 =  ||a-z_1||^2 - ||a_2||^2 < r. \]
Hence $a_1\in E$. This completes the proof of the lemma. 
\end{proof}

\subsection{Relative Albanese reductions}

The goal of this subsection is to prove the following proposition. 
It is a consequence of the relative Albanese reduction, see \cite{Fujiki1983} and \cite{Campana1985}.
We will apply the proposition to the morphism $f'$ obtained in Lemma \ref{lemma:choose-2-form} for the proof of Theorem \ref{thm:uniruledness}.

\begin{prop} 
\label{prop:reduction-Alb}
Let $f\colon X\to Y$ be a fibration between compact complex analytic varieties, 
such that $\dim X>\dim Y$ and that $X$ is a K\"ahler manifold. 
Then, up to blowing up $X$, 
there is a factorization of $f$ into fibrations $p\colon X\to W$ and  $q\colon W\to Y$, 
such that the following properties hold. 
\begin{enumerate} 
    \item $W$ is a compact K\"ahler manifold. 
    \item There is a Zariski open subset $W^\circ \subseteq W$, whose complement is a simple normal crossing divisor, such that $p$ is smooth over $W^\circ$. 
    \item Let $F$ be a general fiber of $p$.  Then $\dim F >0$. 
     In addition, either $H^1(F,\mathcal{O}_F) = \{0\}$ or $F$ has maximal Albanese dimension. 
    \item General fibers of $q$ are non uniruled. 
    In particular, $W$ is not uniruled if $Y$ is not.     
\end{enumerate}  
\end{prop}  

We recall that a compact complex manifold $F$ is said to have maximal Albanese dimension, 
if the Albanese morphism $F\to \mathrm{Alb}(F)$ is generically finite over its image.

\begin{proof}
We first note that, if $X\to W$ is a fibration, then $W$   belongs to the Fujiki's class, 
and hence is bimeromorphic to a compact K\"ahler manifold, see 
\cite[Section IV.3]{Varouchas1989}. 
Thus,  once we have constructed a complex analytic variety $W$ such that the properties (3) and (4) hold, 
we can   blow up $X$ and $W$  further in order  that the properties (1) and (2) are satisfied.  

Let $y\in Y$ be a general point. 
If the  fiber $X_y$ of $f$ satisfies the property  (2), then we just let  $W=Y$.  
Otherwise, we have $H^1(X_y,\mathcal{O}_{X_y}) \neq  \{0\}$.  
Let  
\[X \overset{\xi}{\dashrightarrow} \mathrm{Alb}(X/Y) \overset{\alpha}{\rightarrow} Y\] 
be the relative Albanese reduction 
as in \cite[Th\'eor\`eme 1]{Campana1985}.  
By construction, the restriction $\xi|_{X_y}\colon X_y \to \mathrm{Alb}(X/Y)_y$ is an Albanese morphism of $X_y$, see  \cite[Lemme 2 and the proof of Th\'eor\`eme 1]{Campana1985}.    
It follows that the image $Z$ of $\xi$ has dimension greater than $Y$. 
Up to blowing up $X$, we may suppose  that $\xi$ is a morphism.  
By assumption, $X_y$ does not have maximal Albanese dimension. 
Hence the natural morphism $X\to Z$ has positive relative dimension.    
We note that general fibers of $Z\to Y$ are contained in  tori. 
Hence they are non uniruled.

Let $X\to V \to Z$ be the Stein factorization. 
Then $\dim Y < \dim V < \dim X$ and general fibers of $V\to Y$ are non uniruled.  
Therefore, if general fibers of $g\colon X\to V$ satisfies the property  (3), then we can let $W=V$ (up to proper bimeromorphic models, as mentioned in the first paragraph). 
Otherwise, we apply the same construction as in  the previous paragraph to $g\colon X\to V$. 
Since the relative dimension decreases strictly for such a construction, 
after finitely many steps, we can obtain such a manifold $W$. 
This completes the proof of the proposition.  
\end{proof}

In the case when general fibers have maximal Albanese dimensions, 
we have the following assertion.  
The other case will be studied in the next subsection.

\begin{lemma}
\label{lemma:max-Alb-fiber}
Let $f\colon X\to Y$ be a fibration between compact K\"ahler manifolds. 
Assume that general fibers $F$ of $f$ have maximal Albanese dimension. 
Then $\omega_{X/Y}$ is pseudoeffective. 
\end{lemma}

\begin{proof}
Since $F$ has maximal Albanese dimension, its Kodaira dimension is at least $0$ by \cite[Lemma 10.1 and Theorem 6.10]{Ueno1975}. 
%By Lemma \ref{lemma:semicontinuity}, 
Hence there is an integer $m>0$ such that $H^0(F,\omega_F^{\otimes m}) \neq \{0\}$ for general fibers $F$.  
Thus,   by \cite[Theorem 3.5]{Cao2017}, 
the line bundle $\omega_{X/Y}^{\otimes m}$ admits a positive singular Hermitian  metric.  
It follows that  $\omega_{X/Y}$ is pseudoeffective.  
\end{proof}

\subsection{Morphisms with  relative Albanese dimension zero}

The objective of this subsection is to prove the following proposition.

\begin{prop}
\label{prop:h1=0-psef}
Let $f\colon X\to Y$ be a fibration between compact K\"ahler manifolds.   
We assume the following conditions. 
\begin{enumerate}
    \item  $H^1(F,\mathcal{O}_F)=\{0\}$ for a general fiber $F$ of $f$.  
    \item  $\omega_F$ is pseudoeffective.  
    \item $f$ is smooth over an open subset of $Y$ whose complement is a snc divisor. 
    \item There is a K\"ahler form $\omega$ on $X$ and a closed holomorphic $2$-form $\tau\in H^0(Y,\Omega_Y^2)$, such that the cohomology class $\{\omega + f^*\tau + f^*\overline{\tau}\}$ belongs to the image of  $H^2(X,\mathbb{Z})$.  
\end{enumerate}
Then $\omega_{X/Y}$ is pseudoeffective.   
\end{prop}

We will use the   notation of the proposition throughout this subsection. 
By the same argument of Lemma \ref{lemma:choose-2-form}, 
for any point $o\in Y$, 
there is a Stein open neighborhood  $U\subseteq Y$  of $o$,  
there is a $f$-ample line bundle $\mathcal{L}$ on  $V:=f^{-1}(U)$, 
and a smooth Hermitian metric $h$ on $\mathcal{L}$, 
such that  $\Theta_h = \omega|_V$, 
where $\Theta_h$ is the curvature form of $h$.  
Since $Y$ is compact, we can cover  it  by a finite family of Stein open subsets $\{U_i\}_{i\in I}$ satisfying the previous properties. 
Let $V_i=f^{-1}(U_i)$ and let $(\mathcal{L}_i, h_i)$ be such a Hermitian line bundle. 
We denote $U_{i,j} = U_i\cap U_j$ and $V_{i,j} = V_i \cap V_j$.

\begin{lemma}
\label{lemma:overlap-L} 
Up to multiplying  $\omega$  by a  sufficiently divisible positive integer, 
for any $i,j \in I$, 
there is an isomorphism $\gamma_{i,j} \colon  \mathcal{L}_i|_{V_{i,j}}  \cong \mathcal{L}_j|_{V_{i,j}}$. 
In addition,   there is a real smooth function  $\eta_{i,j}$ on $U_{i,j}$, 
such that $h_i = \gamma_{i,j}^*h_j \cdot e^{-2\eta_{i,j}\circ f}$.  
\end{lemma}

\begin{proof}
Since the index set $I$ is finite, it is sufficient to prove the statement for some fixed pair $i\neq j\in I$.  
The condition (3) implies that $R^1f_*\mathcal{O}_X$ is locally free, 
by the local freeness theorem of \cite{Takegoshi1995} and the relative duality of \cite{RamisRugetVerdier1971} (see also \cite[Theorem 2.6]{Kollar1986II} for the algebraic case).
Hence  $R^1f_*\mathcal{O}_X = 0$ by the condition  (1).  
Since $U_{i,j}=U_i\cap U_j$ is the intersection of   Stein open subsets, it is also Stein. 
Thus $H^1(V_{i,j},\mathcal{O}_{V_{i,j}}) = \{0\}$ after the Leray spectral sequence.   
We note that 
\[c_1(\mathcal{L}_i|_{V_{i,j}}) = c_1(\mathcal{L}_j|_{V_{i,j}}) = \{\omega|_{V_{i,j}}\}  \in H^2(V_{i,j},\mathbb{R}).  
\]
By Lemma \ref{lemma:torsion-line-bundle}, there is some integer $m>0$ 
such that  $\mathcal{L}_i^{\otimes m}|_{V_{i,j}}  \cong \mathcal{L}_j^{\otimes m}|_{V_{i,j}}$.  

Replacing $\omega$ by $ m \omega$, we may assume that $m=1$.   
Then there is a smooth real function $\psi$ on $V_{i,j}$ such that $h_i= \gamma_{i,j}^*h_j \cdot e^{-2\psi}$ on $V_{i,j}$. 
Since   $h_i$ and $h_j$ have the same   curvature form on $V_{i,j}$, 
we deduce that $\psi$ is a pluriharmonic function on $V_{i,j}$.  
Since the  fibers of $f$  are compact, 
by the maximum principle, 
we obtain that  $\psi$ is constant on every fiber. 
Hence there is a continuous function $\eta_{i,j}$ on $U_{i,j}$ 
such that $\psi= \eta_{i,j} \circ f$ on $V_{i,j}$.  
We remark that $\eta_{i,j}$ is pluriharmonic as well.
Hence it is smooth. 
This completes the proof of the lemma. 
\end{proof}

Now we can deduce Proposition \ref{prop:h1=0-psef}.

\begin{proof}[{Proof of   Proposition \ref{prop:h1=0-psef}}]
Let $a>0$ be an arbitrary integer.  
For each $i\in I$,  
since $\omega_F$ is pseudoeffective for general fibers $F$ of $f$, 
%by Lemma \ref{lemma:semicontinuity}, 
there is an integer $b>0$ such that  
\[
H^0(F,  \omega_F^{\otimes ab} \otimes  \mathcal{L}_i ^{\otimes b} ) 
= H^0(F, (\omega_F^{\otimes a} \otimes  \mathcal{L}_i)^{\otimes b} ) \neq \{0\}\]
for a general fiber $F$ over $U_i$. 
%\[f_*\left( (\omega_{X/Y}{^\otimes} \otimes  \mathcal{L_i})^{\otimes b} \right) \neq 0\]
Since $I$ is finite, we can pick the integer $b$ so that this  property holds for all $i\in I$. 
We endow $\omega_{X/Y}^{\otimes ab}|_{V_i}$ the metric $g_i$ constructed in Lemma \ref{lemma:Bergman-indep}, 
with respect to the Hermitian line bundle $(\mathcal{L}^{\otimes b}, h_i^{\otimes b})$. 
Then its curvature current satisfies 
\[
\Theta_{g_i} \ge -\Theta_{h_i^{\otimes b}} = -b \omega|_{V_i}. 
\]
After Lemma \ref{lemma:overlap-L}, 
the second part of Lemma \ref{lemma:Bergman-indep}   implies that the metrics $\{g_i\}_{i\in I}$ are compatible along the overlaps $V_{i,j}$. 
It follows that there is a singular Hermitian metric on $\omega_{X/Y}^{\otimes ab}$, 
such that its restriction on $V_i$ is equal to $g_i$ and that its curvature current satisfies 
\[
\Theta_{g} \ge  -b \omega. 
\]
This shows that 
\[
ab \cdot c_1(\omega_{X/Y}) + b\cdot \{ \omega\} \in H^{1,1}(X,\mathbb{R})
\]
is a pseudoeffective class. 
Since $a$ can be arbitrarily large, 
we deduce that $\omega_{X/Y}$ is pseudoeffective.  
This completes the proof of the proposition. 
\end{proof}

\subsection{End of the proof}

We will complete the proof of  Theorem \ref{thm:uniruledness}. 
For the notion of rational quotient or MRC fibration, 
we refer to   \cite{Campana1992} or \cite{KMM92b}.  
In the setting of compact K\"ahler manifold, see   \cite[Theorem 2.6]{Campana2004ap}.

\begin{lemma}
\label{lemma:mrc-base}
Let $X$ be a uniruled compact K\"ahler manifold and let $f\colon X\dashrightarrow Y$ be the rational quotient or MRC fibration. 
Then $Y$ is not uniruled. 
\end{lemma}

\begin{proof}
It is enough to show that every fibration $g\colon Z\to C$ with rationally connected general fibers  from a compact K\"ahler manifold to a smooth curve admits a section.  
Since a  general fiber $F$ of $g$ is rationally connected, 
it does not have non zero holomorphic differential forms. 
Since in addition $\dim C=1$, we deduce that $H^0(Z,\Omega_Z^2) = \{0\}$. 
Hence $Z$   is projective by Kodaira's embedding theorem.  
Now we can deduce that $g$ admits a section  by \cite[Theorem 1.1]{GraberHarrisStarr2001}. 
This completes the proof of the lemma.
\end{proof}

\begin{proof}[{Proof of Theorem \ref{thm:uniruledness}}] 
Let $\dim X = n$. 
If $X$ is uniruled, then it is cover by free rational curves. 
If $C$ is one of them, 
then the class $\{[C]\} \in H^{n-1,n-1}(X,\mathbb{R})$ of the integration over $C$ is a movable class, 
and we have  $c_1(\omega_X) \cdot \{[C]\} <0$. 
Hence $\omega_X$ is not pseudoeffective. 

Assume from now on that $\omega_X$ is not pseudoeffective.  
Thus, there is some movable class $\alpha\in H^{n-1,n-1}(X,\mathbb{R})$ such that $\mu_\alpha(T_X) >0$. 
We can assume by induction that the theorem holds in smaller dimensions, for it is true if $\dim X =1$.  
If $H^{2,0}(X,\mathbb{C}) = \{0\}$, then $X$ is projective by Kodaira's embedding theorem, 
and it is uniruled by \cite[Theorem 2.6]{BoucksomDemaillyPuaunPeternell2013}. 
Thus we may suppose that $X$ is not projective, and hence  $H^{2,0}(X,\mathbb{C}) \neq  \{0\}$.  

Assume by contradiction that $X$ is not uniruled.   
We note that if $X'\to X$ is a projective bimeromorphic morphism with $X'$ smooth, 
then $\omega_{X'}$ is not pseudoeffective neither.  
Thus, by  Lemma \ref{lemma:kernel-foliation}  and  Lemma \ref{lemma:choose-2-form},  
up to blowing up $X$, 
we may assume that there is a fibration $f\colon X\to Y$ with $0< \dim Y < \dim X$,  
and there is some K\"ahler form $\omega$ on $X$  satisfying the properties in Lemma \ref{lemma:choose-2-form}.  
Since $X$ is not projective, $\{\omega\}$ is not a rational class. 
Thus  $  H^{0}(Y,\Omega_Y^2)  \neq  \{0\}$.

Up to blowing up $X$ further (see Remark \ref{rmk:bimeromorphic}), we can assume that $f$ factors through a morphism  
$p\colon X\to W$  as in Proposition \ref{prop:reduction-Alb}.   
Then general fibers of $p$ are non uniruled. 
Hence by induction hypothesis, their canonical line bundles are pseudoeffective.  
%We note that $p^*\colon  H^{2,0}(W,\mathbb{C}) \to H^{2,0}(X,\mathbb{C})$ is isomorphic as well. 
By applying Proposition \ref{prop:h1=0-psef} or Lemma \ref{lemma:max-Alb-fiber}, 
we deduce that $\omega_{X/W}$ is pseudoeffective. 
Since $\omega_X$ is not pseudoeffective, we obtain that $\omega_W$ is not pseudoeffective neither. 
By induction hypothesis, we deduce that $W$ is uniruled. 
Hence so is $Y$ by the property (4) of Proposition \ref{prop:reduction-Alb}. 

We will now show that $Y$ is rationally connected. 
Assume by contradiction that it is not. 
Let $g\colon Y\dashrightarrow V$ be the rational quotient or the MRC fibration of $Y$. 
Then $V$ is  non uniruled by Lemma \ref{lemma:mrc-base}.  
We may assume that $V$ is  a compact  K\"ahler manifold, see  \cite[Section IV.3]{Varouchas1989}.  
Up to blowing up $X$ and $Y$, we may assume that $g$ and  the natural meromorphic map 
$\overline{f}\colon X\to V$ are morphisms. 
Since general fibers of $g$ are rationally connected, which do not carry non zero holomorphic differential form,  
we see that $g$ induces an isomorphism  $ H^{0}(V,\Omega_V^2) \cong   H^{0}(Y,\Omega_Y^2)$.  
Hence the K\"ahler form $\omega$ on $X$ satisfies the properties of Lemma \ref{lemma:choose-2-form} with respect to $\overline{f}$.  
By applying the same  argument  as the one for $f\colon X\to Y$ to the morphism $\overline{f}$,   
we deduce that $V$ is uniruled. 
This is a contradiction. 

In conclusion, $Y$ is rationally connected.   Thus 
$  H^{0}(Y,\Omega_Y^2)  = \{0\}$, 
which is a contradiction to the third paragraph. 
Hence $X$ is uniruled. 
\end{proof}

\bibliographystyle{alpha}
\bibliography{reference}

\end{document}